\input amstex
\documentstyle{amsppt}
\NoBlackBoxes

\def\O{{\Cal O}}

\def\E{{\Cal E}}

\def\I{{\Cal I}}
\def\J{{\Cal J}}

\def\P{{\Bbb P}}
\def\Q{{\Bbb Q}}
\def\Z{{\Bbb Z}}
\def\C{{\Bbb C}}
\def\F{{\Bbb F}}
\def\S{{\Bbb S}}
\def\R{{\Bbb R}}

\def\f{\varphi}
\def\ra{\rightarrow}
\def\raa{\longrightarrow}
\def\iso{\simeq}

\def\har#1{\smash{\mathop{\hbox to .8 cm{\rightarrowfill}}
\limits^{\scriptstyle#1}_{}}}

\topmatter
\title A view on contractions of higher dimensional
varieties.\endtitle
\author Marco Andreatta and Jaros{\l}aw A.
Wi\'sniewski\endauthor
\leftheadtext{MARCO ANDREATTA and JAROS{{\L}}AW A. WI\'SNIEWSKI}%

%\dedicatory This paper is dedicated to our authors.\enddedicatory

\address  Dipartimento di Matematica, Universit\'a di Trento,
38050 Povo (TN), Italy \endaddress

\email andreatt\@science.unitn.it\endemail

\address Instytut Matematyki UW, Banacha 2, 02-097 Warszawa, Poland
\endaddress

\email jarekw\@mimuw.edu.pl\endemail

%  The following items provide publication information for the AMS-P logo
%\cvol{00}
%\cvolyear{0000}
%\cyear{0000}

%  Math Subject Classifications
%\subjclass Primary: 14E30, 14J40
%Secondary: 14J30, 32H02 \endsubjclass

%\abstract  \endabstract

%  \thanks will become a 1st page footnote.
%  Use \endgraf to indicate a new paragraph; a blank line or \par will
%  be recognized as an error.
\thanks Our contribution to this field was presented during AMS Summer Research
Institute on Algebraic Geometry (July 1995, Santa Cruz). We would like to
thank the organizers of this conference for their great job which includes the
preparation of the present volume. The final version
of this paper was much influenced by referees' remarks which we appreciate
a lot. During the preparation of the paper we
enjoyed hospitality and support of Max-Planck-Institut f\"ur
Mathematik in Bonn.  The first named author was moreover partially
supported by MURST and GNSAGA.  The second named author was partially
supported by an Italian CNR grant (GNSAGA) and by a Polish grant (KBN
GR564 contract 2P03A01208).
We would like to thank all the mentioned above institutions.
\endthanks

\endtopmatter

\document

\subhead Introduction \endsubhead
In this paper we discuss some recent results about maps of
complex algebraic varieties.  A {\sl contraction} $\f:X\ra Z$ is a
proper surjective map of normal varieties with connected fibers.
We will always assume that $X$, if not smooth, has mild singularities;
in particular we assume that a rational multiple of the canonical sheaf,
$rK_X$, is a Cartier divisor. In this hypothesis
the contraction $\f$ is called {\sl Fano-Mori}, or {\sl extremal},
or just {\sl good} if the anticanonical divisor
of $X$, denoted by $-K_X$, is $\f$-ample.
\par

Fano-Mori contractions occur naturally in classification theories of
algebraic varieties starting with the theory of surfaces. This
becomes apparent from the following list of contractions of
smooth surfaces. Namely, suppose that $\f:X\ra Z$ is an extremal
contraction of a smooth surface $X$, then one of the following
occurs:
\item{(a)} $Z$ is a point and $X$ is a  del Pezzo surface;
\item{(b)} $Z$ is a smooth curve and $\f: X\ra Z$ is a conic
or $\P^1$-bundle, in particular every fiber of $\f$ is reduced and
isomorphic to $\P^1$ or to a union of two $\P^1$'s meeting
transversally;
\item{(c)} $Z$ is  a smooth surface (thus $\f$ is birational) and
the exceptional locus consists of disjoint smooth rational curves
with normal bundle $\O(-1)$, thus $\f$ is a composition of blow-downs
of disjoint rational curves to smooth points on $Z$.
\par\noindent
This is a classical result due to the Italian school of algebraic geometry;
for a modern presentation see for instance \cite{B-P-V}.

\noindent
As it follows from the 2-dimensional example, contractions can be either
{\sl birational} or of {\sl fiber type} if $dimZ<dimX$. A birational
contraction $\f$ is called {\sl crepant} if $K_X=\f^*K_Z$. In case of
surfaces crepant contractions lead to {\sl Du Val singularities} which are also
called {\sl simple rational double points} or {\sl A--D--E}-singularities (see
\cite{B-P-V}).
\par

A classical application of Fano-Mori contractions is the general
adjunction theory as presented in \cite{B-S3}.
If $L$ is an ample line bundle on $X$ and
$K_X$ is not nef then one considers an adjoint divisor $K_X+\tau L$,
where $\tau > 0$ is a rational number such that $K_X+\tau L$ is nef
but not ample (rationality of $\tau$ is a theorem of Kawamata, see \cite{K-M-M}
(4.1.1);
such $\tau$ is called {\sl nef} or {\sl threshold} value of the pair
$(X,L)$; see (2.9) and the following).
Then some multiple of $K_X+\tau L$ is spanned (see \cite{K-M-M} (3.1.1)) and the
adjunction mapping of $X$ is the connected part of the Stein
factorization of the map associated to $|m(K_X+\tau L)|$.  The
adjunction mapping is an example of a good contraction.
\par

However, the most natural modern application of extremal contractions
is in the
Minimal Model Program, where the point is to use good contractions and
birational transformations related to them to produce a minimal model
of a given variety. We do not intend to describe how the program works
--- for this we refer the reader to \cite{K-M-M}. We merely note that
the  complete description of (elementary) good
contractions of smooth 3-folds was given by S.~Mori \cite{Mo2} as the base
of the Minimal Model Program in dimension 3  and subsequently, the objective of
the Program in dimension 3, that is the existence of minimal models,
was achieved by studying contractions of varieties admitting
good (terminal) singularities, see \cite{Mo3}.
Contractions of toric varieties were studied for the purpose of the
Program (in arbitrary dimension) by M.~Reid, see \cite{Re1}. Further on,
the theory
of good contractions of terminal threefolds was extended in \cite{Ko-Mo}
and \cite{Ka4}.
Below we recall the main classification result of \cite{Mo2}.

\proclaim {Theorem \rm (Mori, \cite{Mo2})}
Let $\f: X\ra Z$ be an extremal contraction of a smooth projective 3-fold.
Suppose moreover that $\f$ is elementary which means that

\noindent
$Pic(X)/\f^*Pic(Z)\iso{\Z}$.
\item{} If $\f$ is of fiber type then one of the following occurs:
\itemitem{(a)} $Z$ is a point and $X$ is a Fano manifold,
\itemitem{(b)} $Z$ is a smooth curve and $\f$ is a fibration of del Pezzo
surfaces,
\itemitem{(c)} $Z$ is a smooth surface and $\f$ is a conic bundle.
\item{} If $\f$ is birational then its exceptional locus
is an irreducible divisor $E$ and one of the following occurs
\itemitem{(d)} $\f(E)$ is a smooth curve, $Z$ is smooth
and $\f$ is a blow-down morphism
\itemitem{(e)} $\f(E)=z$ is a point, $E\iso\P^2$ and $\O_E(E)=\O(-1)$ or
$\O(-2)$ and then $z$ is a smooth point or a quotient singularity of type
$\C^3/\Z_2$, respectively,
\itemitem{(f)} $\f(E)=z$ and $E$ is an irreducible quadric,
smooth or singular, and then singularity of $Z$ at $z$ is biholomorphic
to the hypersurface singularity of type $x^2+y^2+z^2+t^k$, where
$k=2$ or 3, respectively.
\par\noindent Moreover in cases (e) and (f) the formal neighborhood of $E$
is determined uniquely
(it depends only on the isomorphism class of $E$ and its normal bundle)
and the map $\f$ is the blow-up of $Z$ along the maximal
ideal at $z$.
\endproclaim

In the present survey we focus on Fano-Mori contractions of varieties of
dimension
$\geq 4$. Our task is to understand the local structure of such a contraction
$\f: X\ra Z$. We will choose a point $z\in Z$ and we will consider the
contraction $\f$ on a neighborhood of the fiber $F:=\f^{-1}(z)$. We will
frequently assume that $Z$ is affine and the neighboring fibers are ``not
worse'' than $F$ (e.g.~their dimension does not exceed that of $F$).
Let us note that in dimensions 2 and 3 the fiber $F$ was either a curve
or a divisor (if not the whole space or a point), while from dimension 4
on, both the dimension and codimension of the fiber may be $\geq 2$ which
allows new phenomena to occur.

Contemplating upon the theorem of Mori one may conclude that
the local description of $\f$, which we may aim at, should
contain the following information:
\item{(1)}the description of the geometric fiber $F$ (we recall that
the geometric fiber is just the reduced structure on $\f^{-1}(z)$),
\item{(2)}the description of the conormal sheaf $\J_F/\J^2_F$,
\item{(3)}the description of the scheme theoretic fiber $\tilde F$ which
is defined by the pull-back $\f^{-1}(m_z)\cdot \O_X$ of the maximal ideal
of $z$,
\item{(4)}the description of the local ring $\O_{Z,z}$ (or of its completion)
and of the formal or analytic neighborhood of $F$ in $X$.
\par\noindent In what follows we will try to show that the points (1) --- (4)
constitute a natural sequence of intermediate steps in understanding
of the structure of
$\f$. In many instances it is also very useful to understand
other aspects of the geometry of $\f$ around $F$; this includes
the behavior of the Hilbert scheme of low degree rational curves contracted
by $\f$ (we do not discuss this point in this survey)
and the ``resolution'' of $\f$ in terms of ``known'' morphisms
(blow-ups, blow-downs, special fibrations etc.).

The present paper is organized as follows.  In the first two
sections we survey on two fundamental features of contractions: the
vanishing theorems and existence of rational curves in fibers of good
contractions. These two properties, which in many instances can be
used interchangeably, provide a lot of information on the
contraction. We explain several consequences of these two properties,
one of them is a relative base-point-freeness theorem (1.11 and 1.12)
which, subsequently, has several other applications, e.g.~it implies
the normality of the components of low dimensional fibers (1.10.iii).
The relative base point free theorem allows also to apply conveniently
the theory of deformation of (possibly reduced) embedded rational curves.

In Section 3 we study how the conormal sheaf of a fiber
determines the structure of the contraction around the fiber. Namely,
we will assume that $X$ is smooth and $F$ is a locally complete
intersection in $X$. In this situation, if the conormal bundle
$N^*_{F/X}=\I_F/\I_F^2$ is nef (or better spanned) then the local
structure of $\f$ has several nice properties. For example, it turns out
that in some
instances the pair $(F,N^*_{F/X})$ determines uniquely the formal neighborhood
of $F$ in $X$ and thus the singularity of $Z$ at $z$ (see 3.7) --- as in
the Mori's  theorem above.

In Section 4 we focus on contractions with fibers of dimension $\leq 2$.
A major part of this section is related to our paper \cite{A-W2}. In
particular we
provide a complete description of a birational contraction of a
4-fold with an isolated 2-dimensional fiber (4.7). In the present survey we try
to explain the proof of this classification theorem which, essentially,
goes along the steps (1)---(4) indicated above.

In the fifth section we
present some structure theorems for Fano-Mori contractions with high nef
value. We give also
some extensions of these results to the singular case and to the case
of varieties polarized by ample vector bundles.

Throughout the whole survey we try to avoid technical details of proofs
(actually in many instances we simply omit the proof, giving only a reference).
In some cases however we want to give the reader the feeling of what are the
problems to be dealt with and how to use tools (vanishings and rational curves)
to overcome them --- we then do it by presenting a proof in a special situation.
The reader will certainly notice that some results presented in Sections 1
and 2
overlap and  differ only by apparently technical assumptions.
This is intentional and it is meant to show that cohomology
(vanishings) and Hilbert scheme
(rational curves) techniques are of a similar nature and the difference between
them seems to lay in the range of applicability.

\medskip

The choice and the presentation of the results in this paper
reflects the authors' experience and understanding of the subject.
(Thus, as such, might be considered tendentious. We want however to
apologize for any omission or inaccuracy which might have been
caused by our ignorance or misunderstanding of our colleagues'
work.) Our interest in good contractions was
originally motivated by adjunction theory and classification of Fano
manifolds: \cite{A1}, \cite{A2}, \cite{A-M}, \cite{A-B-W1}, \cite{A-B-W2},
\cite{B-W}, \cite{Wi1}, \cite{Wi2}, \cite{Wi3}.
Our interest
in the local structure of contractions grew together with
understanding of examples and tools, i.e.~vanishings and deformation
of rational curves.  Then, applying Kawamata's technique, which we
learned from \cite{Ka2} and \cite{Ko2}, we obtained in \cite{A-W1} a relative
base-point-freeness theorem which allows to apply a variety of
projective geometry tools to deal with the local structure of a
contraction.  The results about 2 dimensional fibers in Section 4 are
from our paper \cite{A-W2}.

Our original task of \cite{A-W2} was to give a classification of
birational contractions of 4-folds with an isolated 2-dimensional
fiber. While dealing with this problem
we were able to prove the normality of
components of a 2-dimensional fiber (Theorem 1.10) and we understood
how to apply deformations of reducible rational curves (Lemma 2.11)
--- which we learned from Koll\'ar's book \cite{Ko3}.
Our approach to the normal sheaf of a fiber was very influenced
by Mori's fundamental papers \cite{Mo2}, \cite{Mo3} and by vector bundle
examples from \cite{S-W1} and \cite{S-W2}.
The first (yet incomplete) version of the classification
theorem (4.7) was presented in June 1994 during Trento conference of
CIRM.

\medskip
Throughout the paper we work over complex numbers
and use standard notation and definitions
compatible with the references which we provide.

\subhead 1. Vanishing theorems and slicing \endsubhead

The cornerstone of the theory of Fano-Mori contractions is
the vanishing theorem due to Y.~Kawamata, E.~Viehweg and
J.~Koll\'ar (see \cite{K-M-M}, section (1-2), or \cite{E-V}, corollary (6.11)):

\proclaim {Vanishing Theorem (1.1)}
Let $X$ be a variety with at most log terminal singularities
and let $\f:X\ra Z$ be a proper morphism
onto a variety $Z$. Then $R^i\f_*\O_X(K_X)=0$ for $i>dimX-dimZ$.
If $L$ is a line bundle on $X$ such that
$-K_X+L$ is $\f$-ample then $R^{i}\f_*L=0$ for $i > 0$.
\endproclaim

Among the numerous applications of the vanishing let us note the following
theorem about the singularities of the target of a contraction, see
\cite{K-M-M},
(5.1.1) and \cite{Ko1}, Cor.7.4.

\proclaim {Proposition (1.2)} Let $\f:X\ra Z$ be a Fano-Mori contraction of
a variety
with log terminal singularities. Then singularities of $Z$ are
rational (in particular they are Cohen-Macaulay).
\endproclaim

If one assumes that the contraction is equidimensional then the
correlation between the singularities of the domain and of the target
seems to be even closer. Namely we have the following result
whose idea is due to Fujita, see \cite{A-B-W2}, (1.4).

\proclaim {Proposition (1.3)}
Let $X\ra Z$ be a Fano-Mori equidimensional contraction from a smooth $n$-fold
$X$. If $Z$ has at most quotient singularities then it is smooth.
\endproclaim

\demo {Proof} We have by hypothesis a finite Galois cover $g: S \ra Z$ from a
smooth variety $S$ which is \`etale outside the singular locus, see
e.g.\cite{Pr}. Let $G$ be the corresponding Galois group.  Let $Y$ be the
normalization of the fiber product of $X$ and $S$ over $Z$.  Then $p:
Y \ra X$ is \`etale in codimension two.  Thus, by purity of the
branch locus, $p$ is \`etale and $Y$ is smooth.  Then the
equidimensional map $Y \ra S$ is flat and its general fiber is a Fano
manifold.  This implies that $\chi (Y_s,\O_{Y_s}) = 1$ for all $s \in
S$. Let $s \in S$ be a point in which $g$ is not \`etale, $z = g(s)$
and let $H$ be the subgroup of $G$ which is the stabilizer of $s$.
The group $H$ acts on the fiber $Y_s$.  Since $\chi (Y_s,O_{Y_s}) =
1$, this implies, because of the general Riemann-Roch theorem (see
for instance \cite{Fl}, example 18.3.9), that $\sharp H=1$.
Therefore $S = Z$ and $Y = X$, and we are done.
\enddemo

Let us recall that a Fano-Mori contraction  $\f$ is an {\sl
elementary} or {\sl extremal ray} contraction if
$rk(PicX/\f^*PicZ)=1$.

\proclaim {Corollary (1.4)}  Let $\f:X \ra Z$ be an elementary contraction
from a smooth $n$-fold $X$ to a normal surface $Z$. Then $Z$ is smooth.
\endproclaim

\demo {Proof} Since $\f$ is elementary and onto a surface then $\f$ is
equidimensional and $Z$ is (algebraically) factorial so that $K_Z$ is
Cartier. (The factoriality of $Z$ is pretty standard: see e.g.~\cite{K-M-M},
Lemma 5-1-5, or the proof of (2.2) in \cite{B-W}.) Since moreover $Z$ has
only rational singularities then they are Du Val singularities, hence
quotient singularities.
\enddemo

This encourages us to propose the following, probably too optimistic,

\proclaim {Conjecture (1.5)} If $\f:X\ra Z$ is an elementary equidimensional
contraction of a smooth variety to a threefold then $Z$ is smooth.
\endproclaim

\medskip

\noindent {\bf (1.6).  The set-up: } Let assume that $\f :X \ra Z$
is a Fano-Mori or crepant
contraction of a variety $X$ with at most log terminal singularities
onto a normal variety $Z$. We choose a point $z\in Z$ and consider
the fiber $F=\f^{-1}(z)$. Since we are interested in a local description
of $\f$ we will assume that $Z$ is affine and for any $z'\in Z$ we
have $dim(\f^{-1}(z'))\leq dimF$. Let $L$ be a $\f$-ample line bundle.
If $\f$ is a good contraction then we choose $L$ and a positive rational
number $r$ such that $K_X+rL$ is trivial on fibers of $\f$; we also
say that $K_X+rL$ is a good supporting divisor for $\f$. Since $Z$
is affine this is equivalent to $K_X+rL=\O_X$, c.f. (2.9). (Note
that given an extremal contraction $\f$ of a Gorenstein
variety we can always choose $L= -K_X$ and
$r=1$.) If $\f$ is crepant we set $r=0$.

Let $\hat F$ be a scheme structure defined on a fiber $\f^{-1}(z)$ by
global functions, that is $\I_{\hat F}=\f^{-1}\J\cdot\O_X$ for some
ideal $\J$ defining a scheme supported at $z$. In particular, if we
take $\J=m_z$ then $\hat F$ is $\tilde F$, the scheme fiber
structure.

\medskip

The vanishing (1.1) implies vanishing results on the fiber;
the proof of the following
proposition can be found in  \cite{Mo2}, 3.20, 3.25.1,
\cite{Fu2}, 11.3, \cite{An}, 2.2, \cite{Y-Z}, Lemma 4, \cite{A-W2}, 1.2.1.

\proclaim {Proposition (1.7) \rm (Vanishing of the highest cohomology) }
In the situation of (1.6) let $F'$ be a subscheme of $X$
whose support is contained in the fiber
$F$ of $\f$, so that $\f(F') = z$.
If either $t>-r$ or $t =-r$ and $dimF>dimX-dimZ$ then
$$H^{dimF}(F',tL_{|F'})=0.$$
\endproclaim

\proclaim {Proposition (1.7.1)} In the situation
of (1.6) let also $X'\in |L|$ be the zero locus of
a non-trivial section of $L$. Then we have
$$H^{dim (F \cap X')}(F' \cap X',tL_{|F'\cap X'})=0$$
if either $t > -r +1$
or $t=- r +1$ and $dim (F \cap X') \geq dimX-dimZ$.
\endproclaim

\demo{Proof} We like to give here a proof of (1.7.1), the
proof of (1.7) is similar.
First note that
$R^{i}{\f}_*(tL_{X'})=0$ for $i > 0$ and $t > -r+1$ or
$t= -r+1$ and $i \geq dimX -dim Z$;
this follows from the exact sequence
$$0 \raa -L \raa \O_{X} \raa \O_{X'} \raa 0$$
tensorized by $tL$ and from the theorem (1.1).

Now let $\I_{F'\cap X'}$ be the ideal of $F'\cap X'$ in $X'$ and consider
the sequence
$$0\raa \I_{F'\cap X'}\otimes tL \raa \O_{X'}\otimes tL \raa \O_{F'\cap
X'}\otimes tL
\raa 0.$$
Take the direct image $R^{\bullet}\f_*$. Then
$R^i\f_*(\I_{F'\cap X'}\otimes tL)_z=0$ for $i> q:= dim F \cap X'$.
The map $R^q\f_*(tL_{X'})\ra H^q(F',tL_{F'})$
is thus surjective and the proposition follows from
what we have observed at the beginning.
\enddemo

\medskip
 We would like to note that the second part of the argument
is sufficient to prove (1.7) and in fact it
implies the surjectivity of the restriction of the top cohomology
of {\sl any} sheaf for {\sl any} contraction.
In case of resolutions of isolated contractible singularities
Steenbrik \cite{Ste}, 2.14, proved the surjectivity of the restriction of
{\sl any} positive cohomology of the structural sheaf.
It would be actually desirable
to have some extra vanishing for the scheme $\hat F$ (which should reflect
its closer relation to $X$).
The following result is in this direction.

\proclaim {Proposition (1.8)} In the situation of (1.6) assume that
$dimF\geq 2$.
If either $t>-r$ or $t =-r$ and $dimF-1 >dimX-dimZ$ then
$$H^{dimF-1}(\hat F,tL_{|\hat F})=0.$$
\endproclaim

\demo{Proof} The sheaf  $\I_{\hat F}=\f^{-1}\J\cdot\O_X$ is globally generated
so that we have a surjective morphism $\O_X^m\ra\I_{\hat F}$ whose kernel
we denote by ${\Cal K}$. Therefore we have an exact
sequence
$$0 \raa  {\Cal K} \otimes tL \raa tL^{\oplus m}
\raa tL \raa \O_{\hat F}\otimes tL\raa 0.$$
Now we take the direct image $R^{\bullet}\f_*$. Since
$R^i\f_*({\Cal K} \otimes tL)_z=0$ for $i>dim F= :q$ the cohomology group
$H^{q-1}(\hat F,tL_{\hat F})$
is bounded by $R^{q-1}\f_*(tL)_z$ and $R^{q}\f_*(tL)_z$,
which are zero by (1.1).
\enddemo

As an immediate application of the above results we obtain

\proclaim {Corollary (1.9)} If
$F$ is a one dimensional fiber of Fano-Mori or crepant contraction then
$H^1(F,\O_F)=0$. Moreover all irreducible components of $F$ are smooth
rational curves and the graph of $F$, with edges representing
its components and vertices representing their incidence points, is
simply connected. If $F'$ is a component of $F$ then
$H^1(F',\J/\J^2)=0$, where $\J$ is the sheaf of ideals of $F'$.
\endproclaim

\demo {Proof} The vanishing of $H^1(F',\O_{F'})$ from (1.7) implies that
$F'\iso\P^1$.
On the other hand we have
$$0\raa \J/\J^2\raa \O_X/\J^2\raa \O_{F'}\raa 0$$
and since $H^1(X,\O_X/\J^2)=0$ and the map
$H^0(X,\O_X/\J^2)\raa H^0(F',\O_{F'})$ is surjective
the vanishing for $\J/\J^2$ follows. The rest of
the corollary is proved similarly.
\enddemo

\bigskip
The next theorem gives a first description of the fiber $F$ when
$dimF$ is small with respect to $r$. We can think of it as a
relative version of the theorem of Kobayshi-Ochiai characterizing the
projective space and the quadric (see \cite{K-O}).

\proclaim {Theorem (1.10)}  Suppose that the situation is as in (1.6). Let
$S$ be any component of a non trivial (geometric) fiber $F =
\f^{-1}(z)$, let $S'$ be its normalization and $L'$ be the pull back
of $L$ to $S'$. Then $dimF\geq r-1$ and $\geq r$ if
$dimF>dimX-dimZ$. Moreover:
\item{(i)} if $dimF<r$ or $dimF=r>dimX-dimZ$ then $F$ is
irreducible, $r$ is an integer and $(F,L)\iso(\P^{dimF},\O(1))$,
\item{(ii)} if $dimF<r+1$ or $dimF=r+1>dimX-dimZ$ then the Fujita $\Delta$-genus
of the pair $(S',L')$ is zero,
\item{(iii)} if in  (ii) the bundle $L$ is spanned
(see (1.11) and (1.12) below) then $S$ is normal.
\endproclaim

The proof of the above theorem will be sketched in the remaining part
of this section (see also (2.1) in \cite{A2}).
The estimate on $dimF$ as well as the description
of the normalization of $S$ in (i) and (ii) is due to Fujita \cite{Fu1},
Thm.~2.2,
see also \cite{Y-Z}.  The idea of its proof goes back to the celebrated
Kobayashi-Ochiai characterization of the projective space and
quadric. Namely, one has to consider the Hilbert polynomial of the
pair $(\hat S,\hat L)$, where $\hat S$ is a desingularization of $S'$
and $\hat L$ is the pull-back of $L$. The polynomial $\chi (t) =
\chi (\hat S,t\hat L)$ is of degree $dimS$ in the variable
$t$. By the Leray spectral sequence one gets the vanishing of the top
cohomology of $t\hat L$ in the range provided by (1.7), see \cite{Fu1},
proof of 2.2.  On the other hand, by Kawamata-Viehweg vanishing,
$H^i(\hat S,t\hat L)=0$ for $t<0$ and $i<dimS$, because $\hat L$ is
big and nef. Thus we can count zeroes of $\chi(t)$ which gives the
estimate on the degree of $\chi(t)$ and therefore on $dimS$.

In the boundary cases, when the number of zeroes of $\chi(t)$ is
close to $dimS$, one gets precise information on the pair $(\hat S,\hat L)$
and thus on $(S',L')$.
\medskip

One of the key ingredients of the proof of Theorem (1.10)
is the following base
point free theorem (it also allows to use part (iii) of 1.10).
This is the main technical result in our paper \cite{A-W1}
and it is in fact, at least from our view point, a key step in all the study of
contractions with small dimensional fibers.

\proclaim {Theorem (1.11)  \rm (Relative base point freeness)}
Let $\f:X\ra Z$ be a Fano-Mori contraction supported by $K_X+rL$ as in
(1.6).  Let $F$ be a fiber of $\f$.  Assume that $ dim F < r+1$ or,
if $\f$ is birational, that $dimF \leq r+1$. Then $Bs\vert L\vert :=
supp(coker(\f^*\f_*L\ra L))$ does not meet $F$.  Moreover, (after we
restrict $\f$ to a neighborhood of $F$, if necessary) there exist a
closed embedding $X\ra\P^N\times Z$ over $Z$ such that
$L\iso\O_{\P^N}(1)$ (we say that $L$ is $\f$-very ample).
\endproclaim

The idea of the proof of (1.11) goes back to the classical argument
on resolving the base locus of a linear system. In the recent times
this technique was mastered by Kawamata, see e.g.~\cite{Ka1} and also
\cite{K-M-M}.  For a detailed proof of the above theorem we refer the
reader to \cite{A-W1} (see also \cite{Ka2} and \cite{A-B-W2} where special
cases of
this result are proved).  We were inspired by the paper of
Koll\'ar \cite{Ko2}, which provided an effective approach of dealing with
base point loci.  The main ingredient of the proof is again the
theorem (1.1) and we also use the Kobayashi-Ochiai method (to get
``non-vanishing'').

We note an immediate corollary of (1.11):

\proclaim {Corollary (1.11.1)} In the situation of (1.6) assume that
$\f:X\ra Y$ is a Fano-Mori contraction and $K_X$ is Cartier. If either
$dimF=1$ or
$dimF=2$ and $\f$ is birational then $-K_X$ is $\f$-very ample
in a neighborhood of $F$.
\endproclaim

The base point freeness in case of $dimX=4$, $dimZ=3$ was proved
by Kachi in \cite{Kac}.

\proclaim {Theorem (1.12)} Let $\f:X\ra Z$ be an elementary contraction
of a smooth projective 4-fold to a 3-fold. Let $F$ be a fiber of $\f$ of
dimension $2$, then $Bs\vert L\vert$ does not meet $F$.
\endproclaim

The above theorems (1.11) and (1.12) may suggest the following generalization:

\proclaim {Conjecture (1.13)} In the situation of (1.11) if
$dimX -dimZ < dimF \leq (r+1)$ then $Bs\vert L\vert$ does not meet $F$.
\endproclaim

Let us note that (1.11) in the birational case and, more generally,
(1.13) are the best possible,
that is, the inequalities involving $dimF$ and $r$ can not be improved. Indeed,
there exists a del Pezzo surface $S$ such that $-K_S$ has non empty
base locus, then its contraction to the point provides an appropriate
example.
\medskip

With the base point free theorem now (1.10) (i) follows immediately,
i.e.~any component $S$ of the fiber is a projective space and $L_S=
\O(1)$.

\medskip

The relative base point freeness is very useful because of the following

\proclaim {Lemma (1.14) \rm (Slicing argument) }
Suppose that $\f: X\ra Z$ is as in (1.6).  Let $X'$ be a general
divisor from the linear system $\vert L \vert$. Then outside of the
base point locus of $\vert L\vert$, the singularities of $X'$ are not
worse than these of $X$ (e.g.~if $X$ is smooth then outside the base
point locus $X'$ is smooth too) and any section of $L$ on $X'$
extends to $X$. Moreover, if we set $\f':=\f_{\vert X'}$ and
$L'=L_{|X'}$ then $K_{X'}+(r-1)L'$ is $\f'$-trivial. If either $r>1$
or $r \geq 1$ and $\f$ is birational then the map $\f'$ is a
contraction, i.e.~it has connected fibers.
\endproclaim

The use of inductive arguments applying slicing is typical in this
theory and was done classically in the surface case by Castelnuovo,
Enriques and others (see \cite{C-E}). The reader can find a nice
description of this
method in the book of Fujita \cite{Fu2} where it is called Apollonius method.
In \cite{A-W1} we call it {\sl horizontal slicing} to distinguish it from
{\sl vertical slicing} which is the choice of a divisor $X''\subset
X$ defined by a function from $\O_Z= \f_*\O_X$.
\medskip

The following is an example of an immediate application of the
slicing argument which provides the irreducibility part in (1.10).

\proclaim {Corollary (1.15)} Let $\f :X \ra Z$ and $F$ be as in (1.6).
Assume also that either $dimF< r$ or $dimX-dimZ<dimF\leq r$.  Then
the fiber is irreducible and $L^r\cdot F=1$.
\endproclaim

\demo {Proof}
Assume, by contradiction, that the fiber has (at least) two
irreducible components intersecting in a subvariety of dimension $t
\leq (r-1)$. By the base point freeness of $L$, we can choose $t + 1$
sections of $L$ intersecting transversally in a variety with log
terminal singularities and meeting the two irreducible components not
in their intersection.  By construction the map $\f$ restricted to this
variety has non connected fibers and this is in contradiction with
(1.14). The statement about the degree of $F$ is obtained similarly.
\enddemo

\medskip

If a component $S$ of a fiber of a birational contraction is of
dimension 1 then one can use the arguments from the proof of
nonvanishing theorem (3.1) in \cite{A-W1}, pp. 745--746, to extend a pencil
of sections
of $L_S$ to $X$. In such a case applying the above slicing and
connectedness argument one gets

\proclaim {Corollary (1.16)} If a fiber $F$ of a good birational contraction
contains a component of dimension 1 and $r\geq 1$
then $F\iso\P^1$ and $L_S\iso \O(1)$.
\endproclaim

Let us discuss the question of normality in (1.10) (iii). Because
of the vanishing and the slicing argument it follows by the next lemma.

\proclaim {Lemma (1.17)} Let $S$ be an irreducible (reduced) projective variety
and let $L$ be an ample and spanned line bundle on $S$. Assume that for {\bf
every} curve $C$ which is the intersection of $dimS-1$ sections of
$L$ we have $H^1(C,\O_{C}) = 0$. Then $S$ is normal, $L$ is
very ample and the pair $(S,L)$ has $\Delta$-genus zero.
\endproclaim

\demo {Proof} The normality is proved in \cite{A-W2}, (4.3): the crucial
point is to notice
that any 1 dimensional scheme-theoretic intersection of $dimS-1=:s-1$
sections of
$L$ has no embedded point and thus it is Cohen-Macaulay. Namely: let $C_{gen}$
be a general intersection of $s-1$ members of $|L|$. Then $C_{gen}$ is reduced
and since $H^1(C_{gen},\O_{C_{gen}})=0$ it is a smooth rational curve. Now
let us choose
a scheme-theoretic intersection $C:=H_1\cdot\dots\cdot H_{s-1}$
which contains a prescribed point $x\in S$. Since $|L|$ is spanned we may
choose it
so that $C$ is 1-dimensional and generically reduced. We have
$\chi(\O_C)=\chi(\O_{C_{gen}})$
and thus, since $H^1(C,\O_C)=0$, we get $h^0(C,\O_C)=1$. This implies that
$C$ is reduced
and thus Cohen-Macaulay --- so is $S$ at $x$. Now the Serre's normality
criterion
implies that $S$ is normal. The rest is in \cite{Fu2}.
\enddemo

\example{Example (1.18)}
Now we want to give an example of a birational contraction of a
smooth variety (of dimension $6$) containing a non normal 3
dimensional fiber (a del Pezzo threefold).  This example, which was
suggested to us by M.~Reid, shows that the normality claim in (1.17)
(or (1.10)) is the best possible. On the other hand, Kawamata in
\cite{Ka4} provides an example of a 2 dimensional non-normal fiber
of a birational Fano-Mori contraction of a non-smooth variety, as in (1.6),
with $r<1$.
\par
Let $Z$ be the affine variety in $\C^7$ given
as the zero set of the polynomial $$p (x_1, x_2, y_1, y_2, y_3, t_1, t_2) =
x_1^2A + x_1x_2B + x_2^2 C + x_1^3 +x_2^3 + y_1^4 + y_2^4 + y_3^4+ t_1R+t_2S$$
where $A,B,C$ are linear and $R$, $S$ are cubic
 forms in $y_1, y_2, y_3$.  Let then $X$ be
the blow-up of $Z$ along the plane $x_1=x_2=y_1=y_2=y_3=0$,
with the projection $\f : X \ra Z$.  The exceptional
divisor $E$ is hypersurface in $\P^4\times\C^2$ given by the
equation  $$x_1^2A + x_1x_2B +
x_2^2 C + x_1^3 +x_2^3+t_1R+t_2S=0,$$
where, by abuse, $x_i$ and $y_j$ are homogenous coordinates
in $\P^4$.
In particular the fiber of $\f$ over the origin,
$F=E\cap \{t_1=t_2=0\}$, is a non normal cubic which is
singular along the plane $x_1 = x_2 = 0$ and it is irreducible
for a general choice of $A,\ B,\ C$.
The smoothness of $X$ is checked by using affine coordinates in
the blow-up of $\C^7$ in which $X$ is a divisor. For example,
over the set where $y_1\ne 0$ the variety $X$ is given by the
polynomial
$p_1= p (x_1y_1, x_2y_1, y_1, y_2y_1, y_3y_1, t_1, t_2)/y_1^3$,
where this time $x_i$, $y_i$ and $t_i$ are coordinates on the affine piece
of the blow-up. Now evaluating partial derivatives on the plane
$x_1=x_2=y_1=t_1=t_2=0$ we get
$$(\partial p_1/\partial y_1,\ \partial p_1/\partial t_1,\ \partial
p_1/\partial t_2)
=(1+y_2^4+y_3^4,\ R(1,y_2,y_3), \ S(1,y_2,y_3))$$
which is never zero for a sufficiently general choice of $R$ and $S$.
\par
According to a recent preprint of H.~Takagi \cite{Ta1}, this situation
(i.e. non normal fiber) can not occur
if $dimX=4$. (We would like to express our thanks to H.~Takagi for pointing out
a mistake in the previous version of this example.)
\endexample

\medskip

Let us collect the information about two-dimensional fibers of good
contractions.

\proclaim {Theorem (1.19)}
Let $\f :X \ra Z$ be a Fano-Mori contraction as in (1.6) with $r\geq 1$.
Assume moreover that all fibers of $\f$ are of dimension $\leq 2$,
$dimX-dimZ\leq 1$
and $L$ is $\f$-spanned (which is the case if e.g.~$\f$ is birational by
(1.11)).
Let $S$ be a two dimensional component of a (reduced)
fiber $F$ of $\f$. Then $S$ is a
normal surface and the pair $(S,L_{|S})$ is in the following list:
\item{(1)} $(\P^2,\O(e))$, with $e = 1,2$,
\item{(2)} a Hirzebruch surface $(\F_k, C_0 + mf)$ with $m \geq k+1$,
$k\geq 0$, where $f$ is a fiber of the ruling of $\F^r$ over $\P^1$ and
$C_0$ is a section of the ruling such that $C^2_0=-k$,
\item{(3)} a rational normal cone $(\S_k,\O_{\S_k}(1))$, that is a cone over
$\P^1\subset \P^k$ embedded via the $k$-th Veronese map with $k\geq 2$.
\par\noindent
If moreover $\f$ is birational then $F$ is a Cohen-Macaulay scheme of pure
dimension
2 and any two meeting components $S_i$, $S_j$ of $F$ meet along a rational
curve $C_{ij}$
such that $C_{ij}\cdot L=1$.
\endproclaim

\demo {Proof} Normality of $S$ follows directly from (1.10). The list of
pairs $(S,L_S)$
is from \cite{Fu2}. The observation that $F$ is Cohen-Macaulay is in
\cite{A-W2}, (3.4.3),
and it is similar to the proof of normality in (1.17).
Since $F$ is Cohen-Macaulay no two components of
it can have a common isolated point. The pure-dimensionality and the fact that
the intersection of two components is a line follow
from the slicing (1.14) and also (1.9).
\enddemo

We would like to note that even in dimension 3 one can produce examples of
Fano-Mori contractions with a two dimensional fiber $F$
which is a surface from the above list.
For example, a weighted
blow-up $\f:X\ra\C^3$ can have a fiber $F\iso\S_k$ with the value
of $r$ equal to $1+1/k$.
In this case, however, $X$ has an isolated singular point of type
$\C^3/\Z_k(-1,1,1)$.
If $X$ is smooth (or, more generally, has only locally complete intersection
singularities whose locus has zero dimensional intersection
with $F$)
then the possible list of components of the fiber $F$ is smaller, see (4.3.2).

\subhead 2. Rational curves \endsubhead

A fundamental feature of a Fano-Mori contraction is the existence of
rational curves in its fibers. Namely, we have the
following existence theorem due to Mori \cite{Mo1}, \cite{Mo2} in the
smooth case
and extended by Kawamata to the log terminal case in \cite{Ka3}.
We recall
that a rational curve is a curve whose normalization is $\P^1$.
(Although we work over $\C$, we would like to note that if $X$ is smooth
then the existence theorem is true also in positive characteristic; this
concerns also the subsequent results obtained via deformation methods.)

\proclaim {Theorem (2.1) \rm (Existence of rational curves)}
Let $\f:X\ra Z$ be a Fano-Mori contraction of a variety with log terminal
singularities.  Then the exceptional locus of $\f$ is covered by
rational curves contracted by $\f$.
\endproclaim

In this section we study deformations of rational curves
following ideas started with the paper of Mori \cite{Mo1}.  We
discuss only some of the results,
concentrating on the case of smooth $X$. We refer the reader to the
new book of Koll\'ar \cite{Ko3} for general results concerning deformation
of curves.  The following result is from \cite{ibid}, II.1.14.

\proclaim {Theorem (2.2)} Let $C$ be a (possibly reducible) connected curve
such that

\noindent
$H^1(C,\O_C)=0$ and assume that $C$ is smoothable (see
\cite{ibid}, II.1.10 for the definition).  Suppose that $f:C\ra X$ is an
immersion of $C$ into a smooth variety $X$. Then any component of the
Hilbert scheme containing $f(C)$ has dimension $-K_X\cdot C + (n-3)$
at least.
\endproclaim

The above result has several different versions. For example, Mori \cite{Mo1}
proved a version of it for maps of rational curves with fixed points.
An important part of the Mori's proof of the existence of rational
curves is a technique of deforming rational curves with a fixed
0-dimensional subscheme (to ``bend'' these curves) in order to
produce rational curves of lower degree with respect to a fixed ample
divisor (to ``break'' them). In short: if a rational curve can be
deformed inside $X$ with two points fixed then it has to break.
Mori's bend-and-break technique was used by Ionescu (see \cite{Io}, 0.4,
and \cite{Wi2}, 1.1) to prove a bound on the dimension of the fiber.

\proclaim {Theorem (2.3)} Let $\f :X \ra Z$ be a Fano-Mori contraction
of a smooth variety $X$. Let $E$ be the exceptional locus of $\f$
(if $\f$ is of fiber type then $E:=X$) and
let $S$ be an irreducible component of a (non trivial) fiber $F$.
Let $l$ = min $\{ -K_X\cdot C$: $C$ is a rational curve in $S\}$.
Then $dim S + dim E \geq dimX + l -1$.
\endproclaim

The idea of the proof of the above theorem is to use a family of
rational curves obtained by deforming a rational curve in $S$ for
which the minimum $-K_X\cdot C$ is assumed.  If $\f$ is an {\sl
elementary} or {\sl extremal ray} contraction
then the integer $l$ defined above is called the
{\sl length} of the contracted ray.  As an immediate consequence of
the above theorem one gets

\proclaim {Corollary (2.4)} Let $\f :X \ra Z$ be a good contraction of
a smooth variety $X$ and assume that there exists a $\f$-ample line
bundle $L$ such that $\f$ is supported by $K_X +rL$ (see (1.6) or (2.9)).
Then $$dim S + dim E \geq dimX + r -1.$$
\endproclaim

The existence of rational curves implies vanishings on the normalization
of the component. This, together with the Kobayashi-Ochiai argument, was the
base for the following result from \cite{A-B-W2}, Lemma 1.1, (see also
\cite{Zh2}).

\proclaim {Proposition (2.5)} In the hypothesis of the Corollary (2.4)
if the equality holds for an irreducible component $S$ then its
normalization, $\rho : S'\ra S$, is $\P^s$, where $s = dim S$, and
$\rho ^* L = \O(1)$.
\endproclaim

\demo {Proof} We can lift ``minimal length'' curves from $X$ up to the
desingularization $\hat S$ of $S'$ to produce a family of rational
curves on $\hat S$ such that two generic point on $\hat S$ can be
joined by a curve $C$ from this family and $C\cdot \hat L=1$, where
$\hat L$ is the pull-back of $L$.  This implies that $K_{\hat S}\cdot
C\leq -s-1$ and therefore $H^s(\hat S, tL)=H^0(\hat S,K_{\hat
S}-tL)=0$ for $t=-s,\dots,-1$ and by Fujita's version of
Kobayashi-Ochiai result, see \cite{Fu1},  $\hat S=S'\iso\P^s$.
\enddemo

The above results are valid also if we admit some mild singularities,
namely if $X$ has locally complete intersection singularities and $S$
is not contained in the singular locus of $X$. This is because we can
still apply (2.2), see \cite{Ko3}.

\smallskip
One can propose the following conjecture extending the result in (2.5).

\proclaim {Conjecture (2.6)}
In the hypothesis of Theorem (2.3), if equality holds
for an irreducible component $S$ then $S = \P^s$.
\endproclaim

A weaker version of the above conjecture is

\proclaim {Conjecture (2.7)}
If $X$ is a smooth projective variety of dimension $n$ such that
$K_X\cdot C\leq -n-1$ for any complete curve $C\subset X$
then $X\iso\P^n$.
\endproclaim

The above conjecture has drawn attention of many people since the
dawn of Mori theory in \cite{Mo1}. (In particular, the second author of the
present
survey was probably the first who publicized his false proof
of this conjecture in a 1987 preprint.) The conjecture is verified
for $n\leq 3$ in \cite{Wi1} where it is also noted that $PicX\iso\Z$.
\medskip

The inequality (2.3) is very useful when comparing different
elementary contractions. Let us note that if $S_1$ and $S_2$ are
components of fibers of two different elementary contractions then
$S_1\cap S_2$ does not contain a curve. Thus exceptional loci of
elementary contractions of rays of length $l_1$ and $l_2$ do not meet
if $l_1+l_2\geq dimX+3$. This simple observation can be refined in
many ways. A typical result is the following relative version of Mukai
conjecture proved in \cite{B-S-W}.

\proclaim {Proposition (2.8)} Assume that $\f:X\ra Z$ is a fiber type
contraction of a smooth projective variety supported by $K_X+rL$
(see below). If
$r\geq n/2+1$ then $\f$ is elementary unless
$X=\P^{n/2}\times\P^{n/2}$, $L=\O(1,1)$ and $Z$ is a point.
\endproclaim

\medskip
Let us also mention another application of length-fiber-locus
inequality.  For this let us first recall the following.

\proclaim {Definition-Theorem (2.9)}
If $X$ is a projective log terminal variety which is not minimal,
that is $K_X$ is not nef, and $L$ is an ample line bundle
on $X$ then we define
$$\tau (X,L) = {\hbox {\rm min}} \{t \in \R: K_X+tL {\hbox {\rm \ is \
nef}}\}.$$
The number $\tau$ is called the nef value (or the threshold value) of the
pair $(X,L)$ and it is positive and rational
(this last statement is Kawamata's Rationality Theorem -- see \cite{K-M-M},
Theorem (4.1.1)).
Moreover, by  Kawamata-Shokurov Base
Point Free Theorem (see \cite{ibid}, Theorem (3.1.1)) a high multiple of
the divisor $K_X+\tau L$ gives a Fano-Mori contraction $\f$.  The
contraction $\f$ is then called the {\sl nef value
morphism} (relative to $(X,L)$) and $K_X+\tau L$ is called the supporting
divisor
of $\f$.
\endproclaim

Let $\pi: {X}\ra{T}$ be a family of smooth projective
varieties (i.e. $X$ is smooth, $T$ is connected and the map $\pi$ is
smooth and projective) with a $\pi$-ample line bundle ${L}$.  For a
$t\in T$ let $X_t:=\pi^{-1}(t)$ be the associated variety in the
family $\pi$ and let $L_t$ be the restriction of $L$ to $X_t$. The
following was proved in \cite{Wi3}.

\proclaim {Theorem (2.10)} Suppose that $K_X$ is not $\pi$-nef.
Then the function $t\mapsto  \tau (X_t,L_t)$
is locally constant.
\endproclaim

The above theorem is especially nice for Fano manifolds. Namely, let
us assume that $-K_X$ is $\pi$-ample and we choose $\pi$-ample line
bundles $L^1, \dots, L^k$ such that for any $t\in T$ their
restrictions $L^1_t,\dots,L^k_t$ form a real basis for
$N^1(X_t)=PicX_t\otimes_{{\Z}}{\R}$. The choice of
the $L_i$'s provides us with an isomorphism $N^1(X_{t_1})\iso N^1(X_{t_2})$
for any points $t_1,\ t_2\in T$. Then the above theorem implies that
the isomorphism preserves the cone in $N^1(X_t)$ which is spanned by
the classes of ample divisors.
\par
The deformation of rational curves can be also used in case of
crepant contractions. For example, Wilson \cite{Wl} used them to describe
deformations of the ample cone of Calabi-Yau 3-folds.

\bigskip

The existence of rational curves and their deformation can be used to
get information about special fibers of contractions. Suppose
that $\f:X\ra Z$ is a Fano-Mori contraction of a smooth variety and $F$ is
an isolated fiber of $\f$ of dimension $\geq 2$. That is, because $Z$
may be assumed affine, all the fibers of $\f$ except $F$ are of
dimension $\leq 1$. Also, because of (1.10) and (2.12) below, all
1-dimensional fibers of $\f$ are of degree 1, or $\leq 2$, with
respect to $-K_X$, if $\f $ is birational or of fiber type,
respectively.
Let now $C\subset F$ be a rational curve or an immersed
image of a smoothable curve of genus 0. If the degree of $C$ with
respect to $-K_X$ is bigger than that of 1-dimensional fibers of
$\f$, then deformations of $C$ must remain inside $F$ which, in view of
(2.2), provides us with the following useful observation:

\proclaim {Lemma (2.11)} In the above situation
$$dim_{[C]} Hilb(F)\geq -K_X\cdot C+(n-3).$$
\endproclaim

Let us explain why this simple observation is useful for understanding
of the behavior of the morphism $\f$ around $F$.
To show how the above lemma works let us consider the following
situation.

\proclaim {Lemma (2.11.1)}
Let $F$ be an isolated fiber of dimension $\geq 2$ of a
birational contraction $\f: X\ra Z$ of a smooth $n$-fold $X$ and
assume that all the neighboring fibers are of dimension $\leq 1$.
Suppose that $-K_X$ is $\f$-very ample and the exceptional locus of
$\f$ is covered by lines (with respect to $-K_X$).
If there exists a nontrivial decomposition $F=F_1\cup F_2$
then $F_1\cap F_2$ does not contain 0-dimensional components.
\endproclaim

\demo{Proof}
Let $x\in F_1\cap F_2$ be an isolated point
of the intersection. Since $X$ is smooth
$dim_xF_1+dim_xF_2\leq n$.  For $i=1,\ 2$ let $C_i\subset F_i$ be a
a line containing $x$.
Then the variety parameterizing deformations
of $C_i$ inside of $F_i$ with $x$ fixed is of dimension $\leq
dim_xF_i-1$. Let us take $C=C_1\cup C_2$.
Then
$$dim_{[C]} Hilb(F)\leq dim_xF_1+dim_xF_2-2\leq n-2$$ and because
$-K_X\cdot C = 2$ we arrive to the contradiction with (2.11).

\medskip

The assumptions of (2.11.1) may seem rather artificial, nevertheless
they fit very well to the proof of the classification theorem which is
sketch in Section 4. Namely, if $dimF=2$ then by (1.11) $-K_X$ is
$\f$-very ample. In this case, however, (2.11.1) can be also obtained
by a slicing argument --- see (1.19).
\enddemo

\example{Example (2.11.2)}  Let us note that the above conclusion of
(2.11.1) is no longer true if we do not
assume that $\f$ is birational. Namely, we have the following toric
geometry example due to W{\l}odarczyk who generalized our 4-dimensional
example (3.6.2) from \cite{A-W2}.

Let $u_1,\dots, u_k,\ v_1,\dots,v_k$ be a basis of a real vector
space $V$ of dimension $2k$. Moreover let $v_0=-u_0=\sum u_i-\sum
v_i$. We consider a fan $\Delta$ in $V$ containing simplicial cone
$\langle u_1,\dots u_k, v_1,\dots v_k\rangle$ and simplicial cones
$$\langle u_0, u_1,\dots u_{i-1},u_{i+1},\dots u_k, v_1,\dots
v_k\rangle\ \hbox{ and }\ \langle u_1,\dots u_k, v_0, v_1,\dots
v_{i-1},v_{i+1} \dots v_k\rangle$$ for $i=1\dots k$, and the faces of
them. The associated toric variety $X(\Delta)$ is smooth and it
admits a good contraction to a cone $Z$ over
$\P^{k-1}\times\P^{k-1}\subset\P^{k^2-1}$. Namely, we can project the
space $V$ to $V'=V/({\R}\cdot v_0)$. Then the above simplicial
cones are collapsed to the cone spanned on the images of $u_i$s and
$v_i$s and its faces. The contraction $\f:X(\Delta)\ra Z$ has a
unique fiber which consists of two copies of $\P^k$ meeting at one
point, all other fibers are equal to $\P^1$.  (Let us note that
similarly one can produce a contraction of a smooth toric $n$-fold,
such that its general fiber is $\P^1$ and a special fiber consists of
$\P^r$ and $\P^{n-r}$ meeting at one point, however if $r\ne n-r$
then such a contraction is not Fano-Mori.)
\endexample

Let us conclude this section with the following extension of (1.16).
\proclaim {Lemma (2.12)} If a fiber $F$ of a Fano-Mori contraction of a smooth
$n$-fold $X$ contains a component of dimension 1
then $F$ is of pure dimension
$1$ and $-K_X\cdot F\leq 2$.
\endproclaim

\demo{Proof} Let $F'$ be a 1-dimensional component of $F$
($F'$ is a rational curve because of 2.1). Then
$dim_{[F']}HilbX\geq -K_X\cdot F'+(n-3)$ and therefore small
deformation of $F'$ sweep out at least a divisor. More precisely:
taking a small analytic neighborhood of $[F']$ in $Hilb$ and the
incidence variety of curves we can produce an analytic subvariety
$E\subset X$ which is proper over $Z$ such that $F\cap E=F'$ and
$dimE\geq n-1$. This implies that all components of $F$ meeting $F'$
are of dimension 1 and by connectedness of $F$ we see that $F$ is of
pure dimension 1. The bound on the degree can be obtained similarly,
c.f.~(4.1) (note that because of (1.11) $-K_X$ is $\f$-very ample so
that one can apply (2.2) to a curve consisting of two components).
\enddemo

\subhead 3.  The normal bundle of the geometric fiber and the
scheme-theoretic fiber \endsubhead

In the present section we will assume that $\f: X\ra Z$ is a Fano-Mori
contraction of a smooth variety.  Suppose that we already know the
geometric (reduced) structure of the fiber $F=\f^{-1}(z)$ of $\f$.
Then our next task is understanding the scheme-theoretic
structure $\tilde F$ which is defined by the ideal $\f^{-1}m_z\cdot\O_X$.
The understanding of the normal bundle of the geometric fiber $F$ can
be considered an intermediate step in this direction.
\medskip
\noindent
(3.0). We recall that the conormal sheaf (bundle) of $F$ in $X$ is the
quotient $N^*_{F/X} = \I_F /\I_F^2$, where $\I_F$ is the ideal of
$F$. If $F$ is a locally complete intersection then $N^*_{F/X}$ is
locally free over $F$ (a vector bundle) and its dual is the normal
bundle $N_{F/X}$.  We have the adjunction formula: $K_F=(K_X)_{|F} +
det N_{F/X}$.  In particular, if $F$ and $(K_X)_{|F}$ are
fixed then the determinant of the normal bundle is given.  If
$F$ has codimension $1$ the normal bundle is thus fixed.
\par
A very useful example, which we present below, is the contraction to
the vertex of a cone.

\example {Example (3.1)} Let $Y$ be a smooth variety and let ${\Cal L}$ be
a semiample line
bundle over $Y$, i.e.~some multiple of ${\Cal L}$ is spanned by global
sections.
Let $X:=Spec_{Y}(\bigoplus_{k\geq 0} k{\Cal L})$ be the
total space of the dual bundle $\Cal L^*$ with the zero section $Y_0
\subset X$.
Consider the collapsing $\f:X\ra Z$ of $Y_0$ to the vertex $z$ of a
cone $Z$. That is, $Z=Spec(\bigoplus_{k\geq 0} H^0(Y,k{\Cal L}))$ and the
map $\f$ is associated to the evaluation of sections of $k{\Cal L}$.  If $Y$
is Fano and $-K_Y-{\Cal L}$ is ample then the contraction $\f$ is
good. The maximal ideal of the vertex is
$m_z=\bigoplus_{k>0}H^0(Y,k{\Cal L})$. Let us note that
$m_z=\f_*\O_X(-Y_0)$.

The scheme-theoretic structure of the fiber $Y_0$ is defined by sections of
bundles $k{\Cal L}$. More precisely, at a point $y\in Y_0$ the scheme-fiber
structure ideal is generated by functions $s_0^k\cdot s_k$ where
$s_0$ is a local generator of the reduced ideal of $Y_0$ in $X$ (zero
section of ${\Cal L}^*$) and $s_k\in H^0(Y,k{\Cal L})$, $k\geq 1$.
Therefore the fiber
structure coincides with the geometric structure at $y$ if and only
if there exists $s_1\in H^0(Y,{\Cal L})$ which does not vanish at $y$ (or
${\Cal L}$ is generated at $y$).  In particular, the multiplicity of the
fiber is 1, or the fiber structure coincides with the geometric
structure at a general point, if and only if the bundle ${\Cal L}$ has a
non-zero section. In such a case the fiber scheme structure has embedded
components at base-points of the line bundle $\Cal L$.  Let us note
that in this situation, the natural property $H^0(\tilde F,\O_{\tilde
F})=\C$ may be no longer true.

Let us also note that the gradation in the ring $\bigoplus_{k\geq
0}H^0(Y,k{\Cal L})$ may not coincide with the gradation of the maximal
ideal $m_z$ of the vertex. In fact these two coincide if and only if
${\Cal L}$ is projectively normal i.e.~the map $S^k(H^0(Y,{\Cal L}))\ra
H^0(Y,k{\Cal L})$ is surjective for all $k\geq 0$.  If ${\Cal L}$ is ample and
projectively normal then $Z$ is an affine
cone over $Y$ embedded in a projective space by $|{\Cal L}|$.  In such
a case $X$ is the blow-up of $Z$ at $z$ and $m/m^2\iso H^0(Y,{\Cal L})$.
\endexample

We have the following criterion relating the contraction to the
vertex and the blow-up of the vertex.

\proclaim {Lemma (3.2) \rm (\cite{EGA}, (8.8.3))} In the above situation
suppose that ${\Cal L}$ is ample and $H^0(Y,{\Cal L})\ne 0$.
Then $X=Proj_Z(\bigoplus_k m_z^k)$ if and only if ${\Cal L}$ is spanned
over $Y$.
\endproclaim

More generally one can consider a rank $r$
semiample bundle $\E$ on a smooth variety $Y$.
That is $\O_{\P(\E)}(1)$ is semiample, or equivalently,
the symmetric power $S^k(\E)$ is generated
by global sections for $k\gg 0$. Then similarly as above we consider
$X=Spec_Y(\bigoplus_{k\geq 0} S^k\E)$, the total space of the dual
bundle $\E^*$
with the zero section $Y_0$, and the collapsing
$\f:X\ra Z=Spec(\bigoplus_{k\geq 0} H^0(Y,S^k\E))$ of $Y_0$
to the special point $z\in Z$.
If $Y$ is Fano and $-K_Y-det\E$ is ample then the contraction
$\f$ is good. Moreover $\f$ is birational if and only if
the top Segre class of $\E$ is positive. If $\E$ is spanned then the
scheme theoretic fiber $\tilde Y_0$ is reduced and $\f$ factors through
the blow up of $Z$ at $z$ (because of the universal
property of the blow-up).

\bigskip
The subsequent result allows to compare the behavior of an arbitrary
Fano-Mori contraction with a contraction to a vertex; the assumption which
is needed is nefness of the conormal of the fiber. Then the map
behaves similarly as the contraction to the cone because the
sections of the conormal of the fiber extend.

\proclaim {Lemma (3.3)}
Let $\f: X\ra Z$ be a Fano-Mori or crepant contraction of a smooth
variety with a fiber $F=\f^{-1}(z)$. Assume that $F$ is locally
complete intersection and that the blow up $\beta:\hat X\ra X$ of $X$
along $F$ has log terminal singularities. By $\hat F$ we denote the
exceptional divisor of the blow-up. Let ${\Cal L}$ be a line bundle on
$X$ such that $-K_X+{\Cal L}$ is $\f$-big and nef.  If the conormal bundle
$N^*_{F/X}$ is nef then:
\par
\item{(a)} The line bundle $\O_{\hat X/X}(1)=-\hat F$ is $\f\circ\beta$-nef.
\item{(b)} Any section of $N^*_{F/X}$ extends to a function
in $\Gamma(X,\O_X)$
vanishing along $F$, that is the natural map
$$
\aligned
\f^!:m_z&\raa H^0(F,N^*_F)=H^0(\hat F,\O_{\hat F}(-\hat F))\\
m_z\ni f&\mapsto (x\mapsto [f\circ \f]_x\in (\I_F/\I_F^2)_x)
\endaligned
$$
is surjective.
\item{(c)} For $i>0$ and $t\geq 0$ we have vanishing
$H^i(F, S^t(N^*_{F/X})\otimes{\Cal L}) = 0 $.
\item{(d)} Some positive multiple
$\O_{\hat X/X}(k)=-k\hat F$ is $\f\circ\beta$-spanned
and it defines a Fano-Mori contraction over $Z$:
$$\hat\f:\hat X\raa\hat Z=
Proj_Z\bigl(\bigoplus_k (\f\circ\beta)_*\O_{\hat X}(-k\hat F)\bigr);$$
the scheme $\hat Z$ is a blow-up of $Z$ along some ideal of a scheme
supported at $z$.
\endproclaim

\demo {Proof} See \cite{A-W2}.
\enddemo

Subsequently we get the following

\proclaim {Corollary (3.4)} In the situation of the previous lemma
the exceptional set $\hat G$ of the blow-up $\hat Z\ra Z$ is
equal to $Proj(\bigoplus_k H^0(F,S^k(N^*_F)))$ and the map
$\hat\f_{\hat F}:\hat F\ra\hat G$ is defined by the
evaluation
$\bigoplus_k H^0(F,S^k(N^*_F))\ra \bigoplus_k S^k(N^*_F)$.
\endproclaim

So nefness of $N^*_{F/X}$ provides very nice properties for
the local description of the contraction. The situation is even
better if we assume that $N^*_{F/X}$ is spanned.

\proclaim {Proposition (3.5)}
Suppose that $\f: X\ra Z$ satisfies all the assumption of the previous
lemma. Then the following properties are equivalent:
\par
\item{(a)} the bundle $N^*_{F/X}$ is generated by global sections on $F$,
\item{(b)} the invertible sheaf
$\O_{\hat X}(-\hat F)$ is generated by global sections
at any point of $\hat F$.
\item{(c)} $\f^{-1}m_z\cdot\O_X=\I_F$ or, equivalently,
the scheme-theoretic fiber structure of $F$ is reduced, i.e.~$\tilde F=F$.
\item{(d)} there exists a Fano-Mori contraction
$\hat\f : \hat X \ra \hat Z=Proj_Z(\bigoplus_k m_z^k)$
onto a blow-up of $Z$ at the maximal ideal of $z$,
and $\f^*(\O_{\hat Z}(1)) =
\O_{\hat X}(1)$.
\endproclaim

\demo{Proof} The implication (b)$\Rightarrow$(a) is obvious.
(a) implies (b) because of the previous lemma, part (b).
Claims (b) and (c) are equivalent because $\beta_*\O_{\hat X}(-\hat F)=
\I_F$ and $\beta^{-1}(\I_F)=\O_{\hat X}(-\hat F)$.
The implication (b)$\Rightarrow$(d) follows by the universal property of
the blow-up, since by (b) $(\f\circ\beta)^{-1}m_z\cdot\O_{\hat X}=\I_{\hat F}$.
The implication (d)$\Rightarrow$(b) is clear since $\O_{\hat Z/Z}(1)$
is spanned over $\hat Z$.
\enddemo

In some situations the knowledge of the normal bundle $N_{F/X}$ allows
to determine the singularity of $Z$ at $z$.
Let us recall that for a local ring $\O_{Z,z}$ with the maximal ideal
$m_z$ one defines the graded $\C$-algebra $gr(\O_{Z,z}):=\bigoplus_k
m_z^k/m_z^{k+1}$.  The knowledge of the ring $gr(\O_{Z,z})$ allows
sometimes to describe the completion ring $\hat\O_{Z,z}$.  Also, we
will say that a spanned vector bundle $\E$ on a projective variety
$Y$ is p.n.-spanned (p.n. stands for projectively normal) if for any
$k> 0$ the natural morphism $S^kH^0(Y,\E)\ra H^0(Y,S^k\E)$ is
surjective. As we noted while discussing the contraction to the
vertex, projective normality allows us to compare gradings of rings
``upstairs'' and ``downstairs''.

\proclaim {Proposition (3.6) \rm (\cite{Mo2}, 3.32) }
Let $\f:X\ra Z$ be a contraction as above.  Suppose moreover that
$N^*_{F/X}$ is p.n.-spanned. Then $\f_*(\I^k_F)=m^k_z$,\ \
$\f^{-1}(m^k_z)\cdot\O_X=\I^k_F$ and there is a natural isomorphism
of graded $\C$-algebras:
$$gr(\O_{Z,z})\iso\bigoplus_kH^0(F,S^k(N^*_{F/X})).$$
\endproclaim

We omit the proof of the above result referring to Mori, \cite{Mo2},
p.164, who proved it in case when $F$ is a divisor, the
generalization is straightforward. We also note a version of a
theorem of Mori \cite{ibid}, 3.33, which is a generalization of
a Grauert-Hironaka-Rossi result, see \cite{H-R}, \cite{Gr}:

\proclaim {Proposition (3.7)}
Suppose that $F$ is a smooth fiber of a Fano-Mori or crepant
contraction $\f:X\ra Z$ and assume that its conormal bundle
$N^*=N^*_{F/X}$ is nef. If $H^1(F,T_F\otimes S^i(N^*))=H^1(F,N\otimes
S^{i}(N^*))=0$ for $i\geq 1$ then the formal neighborhood
of $F$ in $X$ is determined uniquely and it the same
as the formal neighborhood of the zero section in
the total space of the bundle $N$.
\endproclaim

Also the following assertion is a straightforward generalization of the
celebrated Castelnuovo contraction criterion for surfaces; its proof
is similar to the one of \cite{Ha}, V.5.7, see also \cite{A-W2}.

\proclaim {Proposition (3.8) \rm (Castelnuovo criterion) }
Let $\f :X \ra Z$ be a projective morphism from a smooth variety $X$
onto a normal variety $Z$ with connected fibers (a contraction). Suppose
that $z \in Z$ is a point of $Z$ and $F = \f^{-1} (z)$ is the
geometric fiber over $z$ which is locally complete intersection in
$X$ with the conormal bundle $N^*_{F/X}$.  Assume that for any
positive integer $k$ we have $H^1( F, S^k(N^*_{F/X})) = 0$ (note that in
view of (3.3) this assumption is fulfilled if $\f$ is Fano-Mori contraction,
$N^*_{F/X}$ is nef
and the blow-up of $X$ at $F$ has log terminal singularities). If for
any $k\geq 1$ it is
$S^k H^0(F,N^*_{F/X})\iso H^0(F ,S^k(N^*_{F/X}))$
then $z$ is a smooth point of $Z$ and $dimZ=dimH^0(F,N^*_{F/X})$.
\endproclaim

\medskip
In order to understand higher dimensional fibers of extremal contractions
we slice them by a divisor. Thus we need some kind of
"ascending property" which would allow to compare the properties of
the fiber on the divisor with these on the ambient variety.

Suppose that $\f: X\ra Z $ is a Fano-Mori contraction of a smooth variety,
${\Cal L}$ is an ample line bundle on $X$ such that $-K_X-{\Cal L}$ is $\f$-big
and nef.  Let $F=\f^{-1}(z)$ be a (geometric) fiber of $\f$. Suppose
that $F$ is locally complete intersection. Let $X'\in |{\Cal L}|$ be a
normal divisor which does not contain any component of $F$.  Then the
restriction of $\f$ to $X'$, call it $\f'$, is a contraction.
Moreover, $\f'$ is Fano-Mori or crepant if $-K_{X'}=-(K_X+{\Cal L})_{X'}$
is $\f'$-ample or $\f'$-trivial, respectively (see (1.14)).
The intersection $F'=X'\cap
F$ is then a fiber of $\f'$.  The regular sequence of local
generators $(g_1,\dots,g_r)$ of the ideal of the fiber $F$ in $X$
descends to a regular sequence in the local ring of $X'$ which
defines a subscheme $F\cdot X'$ supported on $F'$, call it
$\bar F'$. We have the following easy

\proclaim {Lemma (3.9)}  The scheme $\bar F'$
is locally complete intersection in $X'$ and $$N^*_{\bar
F'/X'}\otimes_{\O_{\bar F'}}\O_{F'}\iso (N^*_{F/X})_{| F'}.$$
\endproclaim

Let us also note that if the divisor $X'$ has multiplicity 1 along
each of the components of $F'$ then, since a locally complete
intersection has no embedded components, we get $\bar F'=F'$. So, for
the remaining part of this section assume that actually $\bar F'=F'$.
Now we discuss the following extension property.  Let us consider a
point $x\in F'$. Suppose that the ideal of $F'$, or equivalently
$N^*_{F'/X'}$ is generated by global functions from $X'$. That is,
there exist global functions $g'_1, \dots g'_r\in \Gamma(X',\O_{X'})$
which define $F'$ at $x$. Then, since $H^1(X,-{\Cal L})=0$ these functions
extend
to $g_1, \dots g_r\in \Gamma(X,\O_{X})$ which define $F$. Thus passing from
the ideal $\I$ to its quotient $\I/\I^2$ we get the first part of

\proclaim {Lemma (3.10)}
If $N^*_{F'/X'}$ is spanned by global functions from $\Gamma(X',\O_{X'})$ at
a point $x\in F'$ then $N^*_{F/X}$ is spanned at $x$ by functions from
$\Gamma(X,\O_{X})$.
If $N^*_{F'/X'}$ is spanned by global functions from $\Gamma(X',\O_{X'})$
everywhere on $F'$ then $N^*_{F/X}$ is nef.
\endproclaim

\demo {Proof} We are only to proof the second claim of the lemma. Since
$F'\subset F$ is an ample section then the set where $N^*_{F/X}$ is not
generated by global sections is finite in $F$. Therefore the restriction
$(N^*_{F/X})_{|C}$ is spanned generically for any curve $C\subset F$ and
consequently it is nef.
\enddemo

If the fiber is of dimension 2 and $\f$ is birational
then we have a better extension property.

\proclaim {Lemma (3.11)} Let $\f:X\ra Z$ be a Fano-Mori
birational contraction of a smooth variety
with a 2-dimensional
fiber $F$ which is a locally complete intersection.
As usually $L=-K_X$ is a $\f$-ample line bundle which can be
assumed $\f$-very ample (see (1.11)).
Then the following conditions are equivalent:
\item{(a)}$N^*_{F/X}$ is generated by global sections at any point of $F$
\item{(b)} for a generic (smooth) divisor $X'\in |L|$ the bundle
$N^*_{F'/X'}$ is generated by global sections at a generic point of
any component of $F'$.
\endproclaim

\demo{Proof} The implication (a)$\Rightarrow$(b) is clear. To prove the converse
we assume the contrary. First let us note that, since $F'$ is a union
of smooth rational curves with no cycle then any bundle on $F'$
spanned at a generic point of any component is spanned everywhere
(the proof of this is an easy induction with respect to the number of
irreducible components in $F'$).  Now let $S$ denote the set of
points on $F$ where $N^*_{F/X}$ is not spanned. Because of the
extension property (3.10) and what we just said, the set does not
contain $F'$ and thus it is finite. Now we choose another smooth
section $X'_1\in |L|$ which meets $F$ along a (reduced) curve $F'_1$
containing a point of $S$.  (We can do it because $L$ is
$\f$-spanned.) The bundle $N^*_{F'_1/X'_1}$ is generated on a generic
point of any component of $F'_1$ so it is generated everywhere but
this, because of the extension property, implies that $N^*_{F/X}$ is
generated at some point of $S$, a contradiction.
\enddemo

\subhead 4. Good contractions with small fibers \endsubhead

In the present section we will apply results of the preceding
sections in order to discuss good contractions with fibers of small
dimension.  Let $\f :X \ra Z$ be a Fano-Mori (or crepant) contraction of a
variety $X$ which throughout this section will be smooth.  We choose
a fiber $F=\f^{-1}(z)$ and we set $L := -K_X$ if $\f$ is Fano-Mori.  Since
we want to describe the contraction in a neighborhood of $F$, the
target $Z$ will be assumed affine, c.f. (1.6).  One should bear in
mind that the results which we state below concern $\f$ over {\sl some}
neighborhood $Z$ of $z$, so it can be shrunk if necessary.

As it was discussed in the introduction, we want to classify pairs
$(F,L_F)$, describe the conormal sheaf of $F$ in $X$, find out the
singularity of $Z$ at $z$ and describe an analytic or formal
neighborhood of $F$ in $X$. We first assume that $F$ is a fiber of
dimension 1; the following theorem is due to S. Mori in dimension $3$
and to T. Ando for $n\geq 4$.

\proclaim {Theorem (4.1)} Let $\f :X \ra Z$ be a Fano-Mori contraction of
a smooth variety $X$ of dimension $n$. Let $F=\f^{-1}(z)$ be
a fixed fiber and assume that $dim F = 1$.
\item{(1)} If $\f$ is birational then $F$ is irreducible, $F \iso \P^1$,
$-K_X\cdot F=1$ and its normal bundle is $N_{F/X}=\O(-1) \oplus\O^{(n-2)}$.
The target $Z$ is smooth and $\f$ is a blow-up of a smooth codimension 2
subvariety of $Z$.
\item{(2)} If $\f$ is of fiber type then
$Z$ is smooth and $\f$ is a flat conic bundle. In particular one of the
following is true:
\item{(i)} $F$ is a smooth $\P^1$ and $-K_X\cdot F=2$, $N_{F/X}\iso\O^{(n-1)}$;
\item{(ii)} $F=C_1\cup C_2$ is a union of two smooth rational curves meeting
at one point and  $-K_X\cdot C_i=1$, $(N_{F/X})_{|C_i}\iso \O^{(n-1)}$,
$N_{C_i/X}\iso \O^{(n-2)}\oplus\O(-1)$ for $i=1,\ 2$;
\item{(iii)} $F$ is a smooth $\P^1$, $-K_X\cdot F=1$ and the fiber
structure $\tilde F$ on $F$ is of multiplicity 2 (a non reduced conic);
the normal bundle of $\tilde F$ is trivial while $N_{F/X}$ is either
$\O(1)\oplus\O(-1)^{(2)}\oplus\O^{(n-4)}$ or
$\O(1)\oplus\O(-2)\oplus\O^{(n-3)}$
depending on whether the discriminant locus of the conic bundle is smooth
at $z$ or not.
\endproclaim

This theorem was proved by Ando (see \cite{An} and also \cite{Mo2}) using the
vanishing theorem, especially (1.7) and (1.9). In  \cite{A-W2} we
give parallel arguments based on the theory of
deformations of rational curves.  We also note that this theorem was
generalized to the case of a variety $X$ with
terminal Gorenstein singularities by Cutkosky (see \cite{Cu}) for $n=3$,
and by Mori and Koll\'ar (see \cite{Ko-Mo}, (4.9) and (4.10.1)) for $n>3$.
The case of an extremal contraction of a 3-fold $X$ with terminal non Gorenstein
singularities is much more difficult; this was discussed in the celebrated paper
of Mori \cite{Mo3} and in \cite{Ko-Mo}.

\medskip
For crepant contractions the result we have is not as satisfactory as
the one for Fano-Mori contractions.  Let us mention a special
case of 3-folds which will be used later.

\proclaim {Proposition (4.2)}
Let $\f: X\ra Z$ be a crepant contraction of a smooth 3-fold $X$ with
a 1-dimensional fiber $F$.  If $F$ is irreducible then $F\iso \P^1$
and the conormal bundle $N^*_{F/X}=\I_{F}/\I^2_{F}$ is isomorphic to
either $\O(1)\oplus\O(1)$ or to $\O\oplus\O(2)$ or to
$\O(-1)\oplus\O(3)$; note that in the first two cases, because of
(3.5), we have $\I_F=\f^{-1}m_z\cdot\O_X$. If $F$ consists of two
components, $F=F_1\cup F_2$, then each of the components is a smooth
$\P^1$, they meet transversally at one point and the conormal
$N^*_{F/X}$ is locally free and restricted to each
component $F_i$ it is isomorphic to either $\O\oplus\O(1)$ or
$\O(-1)\oplus\O(2)$.  If the restriction to both components is
$\O\oplus\O(1)$ then $\I_{F}=\f^{-1}m_z\cdot \O_X$.
\endproclaim

\demo{Proof} The first part is now well known, see \cite{La} and also
\cite{C-K-M} (16.6),
for the second part see \cite{A-W2}, (5.6.3).
\enddemo

\bigskip
Now we turn to the case of a 2-dimensional fiber.
We follow the arguments of our paper \cite{A-W2}.
The case of small contractions
of 4 folds was resolved by Kawamata in \cite{Ka2}. In this paper Kawamata
proved a
relative base point free theorem (for this case) and studied the conormal
bundle of the fiber; as the final result he obtained a complete description
of this case, including the description of its flip. We were strongly inspired
by this paper.

\bigskip
\noindent
(4.3). Let $\f :X \ra Z$ be a Fano-Mori contraction of a smooth variety of
dimension $n$ with all fibers of dimension less or equal than $2$
and moreover $dimX-dimZ\leq 1$.
Because of (2.12) any two dimensional fiber $F$ is of pure dimension 2. Let
us set $L := -K_X$. Then $L$ is $\f$-spanned  in the
birational case and in the fiber type case if $n=4$
and $\f$ is elementary, see (1.11) and
(1.12), respectively.  In the remaining cases we have to assume that
$L$ is spanned (see Conjecture 1.13).  The list
of possible components $S$ of $F$ is provided in (1.19).

Assume now that all fibers except $F$ have dimension $ \leq 1$. We
will say, for brevity, that $F$ is an isolated two dimensional fiber.
The general non trivial fiber of $\f$ is a line or a conic (relative to
$L$) in the
birational or fiber type case, respectively (see (4.1)).

The list (1.19) of possible components of $F$ is
redundant in this hypothesis as it can be shown by
using deformations of rational curves and Lemma (2.11).
Choosing properly the curve $C$ in the lemma we
can rule out many of possible pairs $(S,L_S)$ which
were listed above. We will give just an example
which explains our argument without discussing numerous cases.
Suppose that $S$ is a component of $F$ and $S\iso\S_r$
where $r\geq 3$. Then as the curve $C$ we take the union of general
$r-1$ lines passing through the vertex of $\S_r$ (the lines are general so
that none of them is contained in any other component of $F$). Then
the deformation argument implies that the curve $C$ has to move out
of $S$ for $n\geq 4$, and thus $\S_r$ can not be a component of $F$,
for $r\geq 3$ if $\f$ is birational and for $r\geq 4$ if $\f$ is of
fiber type.

Also this way for a reducible fiber $F$ we can
limit the possible combination of irreducible components of $F$
--- a typical argument is presented in the proof of (2.11.1).
We note that the very ampleness of $L = -K_X$ as well as
the precise description of the components of the fiber allow us
to choose properly the curve which satisfies
the assumptions in (2.2).

\bigskip
With such arguments we prove the following

\medskip
\noindent
\proclaim {Proposition (4.3.2) \rm (\cite{A-W2}, Sect.~4)}
Let $\f: X\ra Z$ be a Fano-Mori contraction
of a smooth $n$-fold $X$ with an isolated 2-dimensional fiber $F$; let
$L=-K_X$. If $\f$ is birational we have the following possibilities
for the pair $(F, L_F)$

{\baselineskip 20 pt
\settabs\+  $(\P^2, \O(1))$ and $n\leq 6$\ \ \
& $(\P^2\cup \P^2, \O(1))$$C$ \ \ \
& $\F_2 \cup \P^2$; $L_{|\P^2} = \O(1)$ \ \ &
\cr
\+ $n\geq 5$&  n = 4 \ \ & n=3
\cr
\medskip
\hrule
\medskip
\+$(\P^2,\O(1))$, $n\leq 6$ &$(\P^2,\O(1))$ & $(\P^2, \O(1)),\ \  (\P^2, \O(2))$
\cr
\+ &$(\F_0,C_0+f)$ & $(\S_2,\O_{\S_2}(1))$
\cr
\+ & $(\S_2,\O(1))$ & $(\F_0 , C_0 + f)$,\ \ $(\F_1, C_0 + 2f)$
\cr
\+ & $(\P^2\cup \P^2, \O(1))$& $\P^2 \cup_{C_0} \F_2$,  $L_{|\P^2} = \O(1)$,
$L_{|\F_2} = C_0 + 3f$
\cr
}

\medskip\noindent If $\f$ is of fiber type and $L$ is $\f$-spanned then
we have the following possibilities for the pair $(F, L_F)$

{\baselineskip 20 pt
\settabs\+  $(\P^2, \O(1))$ and $n\leq 6$\ \ \ \ &  n = 4, F irreducible \
&  n=4,  F reducible & $\F_2 \cup \P^2$; $L_{|\P^2} = \O(1)$
\cr
\+ $n\geq 5$ &$n = 4$, irreducible\ \ \ \   &$n=4$, reducible &\ \  n=3
\cr
\medskip
\hrule
\medskip

\+ $(\P^2,\O(1))$, $n\leq 7$ & $(\P^2,\O(1))$ & $\P^2\cup\P^2$ & $(\F_0,
C_0 + 2f)$
\cr
\+ $(\F_0,C_0+f)$, $n=5$ & $(\P^2,\O(2))$ & $\P^2\bullet\P^2$ &
$\F_0 \cup \F_1,\  L_{\F_0}=C_0+f$
\cr
\+ $(\P^2\cup \P^2, \O(1))$, $n=5$ & $(\S_2,\O(1))$ &$\P^2\cup\F_0$
$\ \ \ \ \ \ \ \ \ \ \ \ \ \ \   L_{\F_1}=C_0+ 2f$\cr
\+  & $(\S_3,\O(1))$ &$\P^2\cup_{C_0}\F_1$ &
\cr
\+  &$(\F_1,C_0+2f)$ & $\P^2\cup\S_2$
\cr
\+  &$(\F_0,C_0+f).$ & $\P^2\cup\P^2\cup\P^2$
\cr
\+  &  & $\P^2\cup_{C_0}\F_0\cup_{f}\P^2$
\cr
}
\medskip\noindent
In the above list the components of reducible fibers have a common line
(in some cases we point out which line is it) with the unique exception
of two $\P^2$s which meet at a point --- we denote this union by $\bullet$.
(We suppress the description of $L$ whenever it is clear.)
\endproclaim

\bigskip
The proof of the proposition is contained with all details in
Section 4 of our paper \cite{A-W2}.
We note that for $n=3$ we obtain Mori's theorem stated in the introduction
(actually we obtain a slightly
more general result since we do not assume that $\f$ is elementary).
Moreover if $n=3$ then by the adjunction formula the normal bundle of
$F$ is completely determined by the pair $(F,-K{_X}_{|F})$ --- thus
applying (3.6), (3.7) and (3.8) we can describe the singularity of
$Z$ in $z= \f(F)$ and we have the unicity of $\f$ in a formal neighborhood
of $F$ by computing the cohomology obstructions.
The case $n=3$ and non-elementary $\f$ of fiber type was studied
in \cite{B-S1} and \cite{Bs}.

Parallel, also Y.~Kachi \cite{Kac} described the geometric structure of
possible isolated 2 dimensional fibers of a contraction from a 4-fold
to a 3-fold. In this case the list of Kachi coincides with ours (except
one case which we are able to exclude). We would like to note that
our argument relies on the assumption that $-K_X$ is $\f$-spanned,
which is proved by Kachi in case of 4-folds  (see (1.12)),
on the other hand
in his proof Kachi uses our results for the birational case.
His strategy of the proof is different from ours'
and it depends on what he calls
Campana-Koll\'ar-Miyaoka-Mori rationally chain
connectedness technique. In short: Kachi considers limits
of general conics (fibers of $\f$)
in the fiber $F$ and proves that except
the case
$F=\P^2\bullet\P^2$ any two points in $F$ can be joint by such a limit conic.
\par
Another very recent preprint of Kachi \cite{Kac1}
deals with the case of a small contraction of a 4-fold with isolated
locally complete intersection singularities. He proves that in this
case the exceptional locus is an isolated $\P^2$ contracted to a point
(which can be proved by arguments similar to these used to prove
4.3.2 --- because Koll\'ar's Theorem 2.2 works in this case too) and
moreover he describes the flip in this case. He provides an
example of such a contraction of a non-smooth variety and compares it
with Reid's explicit construction of a flop \cite{Re}.
\par
Also, we would like to mention a preprint of H.~Takagi \cite{Ta}
in which the case of a birational
elementary contraction of a smooth 4-fold with non-isolated two dimensional
fiber is treated.

\remark {Remark (4.4)}
Let us note that for almost all the above possibilities we can
construct examples with appropriate isolated 2-dimensional fiber, see
Section 3 of \cite{A-W2} and the following (4.8).  However, there are some
exceptions
for which we were unable to construct examples and we do not expect
that all of them exist.  This concerns only fiber type contractions
and reducible
fibers: $\P^2\cup\P^2$ for $n =5$ and  $\P^2\cup\S_2$,
$\P^2\cup\P^2\cup\P^2$ and $\P^2\cup\F_0\cup\P^2$ for $n=4$.
\endremark

\bigskip

>From now on we concentrate on the case of birational contractions of
smooth 4-folds. We note first that the fiber $F$ is then locally
complete intersection in $X$ and its blow-up has terminal
singularities and thus we can apply results of the
previous section. We want to determine the normal bundle and the
scheme fiber structure of $F$; this is not as straightforward as
in the case $n=3$ since $codim_X F =2$. We will apply slicing technique. Let
$X'\in |L|$ be a normal divisor which does not contain any component
of $F$.  Then the restriction of $\f$ to $X'$, call it $\f'$, is a
crepant contraction and the intersection $F'=X'\cap F$ is then a
fiber of $\f'$. Now, applying previous results, namely
the information about crepant contractions (4.2) and the ascending
property (3.11), we get the following

\proclaim {Lemma (4.5)}
Let $\f:X\ra Z$ be a Fano-Mori birational contraction of a smooth 4-fold with an
isolated 2-dimensional fiber $F=\f^{-1}(z)$.
Then the fiber structure $\tilde F$ coincides with the geometric structure $F$
unless one of the following occurs:
\item{(a)} the fiber $F$ is irreducible and the restriction of
$N_{F/X}$ to any smooth curve $C\in |L_{|F}|$ is isomorphic to
$\O(-3)\oplus\O(1)$,
\item{(b)} $F=\P^2\cup\P^2$ and the restriction of $N_F$ to any line
in one of the components is isomorphic to $\O(-2)\oplus\O(1)$.
\endproclaim

We can now give the description of the normal bundle of the
exceptional fiber. First, however, let us recall that a rank 2 spinor bundle
can be defined on a 2 dimensional quadric as follows. Let ${\Cal S}$
be the universal sub-bundle on the smooth quadric $\Q^4$ which is the
same as the Grassman manifold of lines in $\P^3$. The restriction of
${\Cal S}$ to any plane section of $\Q^4$ will be called
${\Cal S}$ again. This way we define ${\Cal S}$ on $\P^1\times\P^1$, on
the quadric cone $\S_2$ and on the reducible quadric $\P^2\cup\P^2$.
We note that our definition of the spinor bundle on the smooth quadric
$\P^1\times\P^1$ differs from the usual one, which is just a line bundle
associated to one of the rulings, while the rank 2 vector bundle we have just
defined is a sum of such line bundles.

\proclaim {Theorem (4.6)}
Let $\f:X\ra Z$ be a Fano-Mori birational contraction of a smooth 4-fold
with an isolated 2-dimensional fiber $F=\f^{-1}(z)$.  Then
the conormal $N^*_{F/X}$ is spanned by global sections and the scheme
fiber
structure $\tilde F$ coincides with the geometric structure $F$.
Moreover if
$F=\P^2$ then $N^*_{F/X}$ is either $\O(1)\oplus\O(1)$ or
$(T\P^2(-1)\oplus\O(1))/\O$, or $\O^{\oplus 4}/\O(-1)^{\oplus2}$.  If $F$ is a
quadric (possibly singular or even reducible) then $N^*_{F/X}$ is
${\Cal S}(1)$, where ${\Cal S}$ is the spinor bundle.
\endproclaim

\demo{Proof} (See \cite{A-W2} for details.) The first part of the theorem
follows from
the preceding lemma after discussing the possible exceptions described
there.  If $F=\P^2$ this is particularly easy: if $F \not= \tilde F$
then for every line $C \subset \P^2$ we would have $(N_{F/X})_{|C} = \O(-3)
\oplus \O(1)$. Then,  because of the well known theorem of Van de Ven,
$N_{F/X}$ would be decomposable and in particular
$h^0(N_{F/X})-h^1(N_{F/X})>0$.  Similarly one obtain the same
inequality for irreducible quadric, namely using results on vector
bundle due to Forster, Hirschowitz and Schneider; see \cite{F-H-S} and
\cite{A-W2}.  This inequality, however, because of the theory of
deformation (see e.g. \cite{Ko3}), would imply that $F$ moves in $X$ which
contradicts our assumption that $F$ is an isolated 2-dimensional
fiber. The case of a reducible quadric requires, apart of Van de Ven's
theorem, some more discussion on blow-ups, see \cite{A-W2}.

For the second part of the theorem we can use the classification of
spanned vector bundle with low first Chern class on Del Pezzo
surfaces, (see \cite{S-W1}, \cite{S-W2}). We check which bundle among
spanned vector
bundles with the appropriate $c_1$ is actually the conormal of an
isolated 2-dimensional fiber.  We use also the following
observations.  Because of the deformation argument we know that
$h^0(N_{F/X})-h^1(N_{F/X})\leq 0$ and thus $N^*_{F/X}$ can not be
decomposable with a trivial factor. On the other hand none of the
bundles with $s_2=c_1^2-c_2=0$ can occur. Indeed, because of
the classification of such bundles, see \cite{ibid}, we know that
$H^0(S^n(N^*_{F/X})) = S^n(H^0(N^*_{F/X}))$ and $dimH^0(N^*_{F/X})=3$,
thus  by (3.8) the
contraction should be to a 3-dimensional smooth point contrary to the
fact that $dimZ=4$.  Let us note that the cases with $s_2 = 0$ do
occur in the fiber type contractions; see examples (3.5.5) in \cite{A-W2}.
The case of the reducible quadric needs, again, some extra care and
we refer the reader to \cite{ibid} for the discussion of this case.
\enddemo

The next theorem is the main classification result in \cite{A-W2}.

\proclaim {Theorem (4.7)} Let $\f:X\ra Z$ be a birational Fano-Mori contraction
from a smooth variety $X$ of dimension $4$ onto a normal variety $Z$.
Let $F=\f^{-1}(z)$ be a (geometric)
fiber of $\f$ such that $\hbox{dim}F=2$. Assume that all other
fibers of $\f$ have dimension $<2$ and all components of the
exceptional locus $E$ of $\f$ meet $F$ (again, this may be achieved by
shrinking $Z$ to an affine neighborhood of $z$ and restricting $\f$
to its inverse image, if necessary).
\par\noindent
If $\f$ is not divisorial then $E=F\iso\P^{2}$
and its normal is $N_{F/X} = \O(-1) \oplus \O(-1)$.
In this situation the {\it flip} of $\f$
exists. (This was studied in \cite{Ka2}).
\par\noindent
If $E$ is a divisor then $Z$ as well as $S:=\f(E)$ is smooth outside
of $z$. Moreover, outside of $F$ the map $\f$ is a simple blow-down
of the divisor $E$ to the surface $S\subset Z$.  The scheme theoretic
fiber structure over $F$ is trivial, that is the ideal $\I_F$ of $F$
is equal to the inverse image of the maximal ideal of $z$. The fiber
$F$ and its conormal $N^*_{F/X}=\I_F/\I_F^2$ as well as the
singularity of $Z$ and $S$ at $z$ can be described as follows

\medskip
\settabs\+$\P^1 \times \P^1$\ \ \ \ \
&$T_{\P^2} (-1) \cup \O \oplus \O(1)$  \ \ \ \ \ \
&quadratic sing. \ \ \ &cubic singularity \ \cr
\+ $F$  &$N^*_{F/X}$ &$Sing Z$ &$Sing S$\cr
\hrule
\smallskip
\+$\P^2$ &$T(-1)\oplus\O(1)/\O $
&cone over $\Q^3$&smooth \cr
\+$\P^2$ &$\O^{\oplus 4}/\O(-1)^{\oplus 2}$
&smooth &cone over rational twisted cubic\cr
%$(xt-yz, x^2y-z^2,y^2-zt)$ \cr
\+ quadric &spinor bundle from $\Q^4$ &smooth&non-normal\cr
\smallskip
\hrule
\smallskip\noindent
The quadric fiber can be singular, even reducible, and here is the
refined description of its conormal bundle and the ideal of $S$
\medskip
\settabs\+ quadric cone \ \ \ \ \ & the normal bundle\ \ \ \ \ \
 \ \ \ \ \ \ \ \ \  &
the ideal of $S$ in $\C[[x,y,z,t]]$\cr
\+ quadric & conormal bundle & the ideal of $S$ in $\C[[x,y,z,t]]$\cr
\smallskip
\hrule
\smallskip
\+$\P^1 \times \P^1$ &$\O(1,0)\oplus\O(0,1)$
&$(xz,xt,yz,yt)$\cr
\+ quadric cone & $0\ra\O\ra N^*\ra \J_{line}\ra 0$
 & generated by 5 cubics \cr
\+$\P^2 \cup \P^2$ &$T_{\P^2} (-1) \cup (\O \oplus \O(1))$
& generated by 6 quartics\cr
\smallskip
\hrule
\smallskip
\noindent
Moreover all the above cases exist.
\endproclaim

For the complete proof we refer the reader to \cite{A-W2}. The description
of all pairs $(F, N_{F/X})$ is in (4.6).  The rest of the theorem is proved
with the help of results of Section 3.

We give now a sketch of the proof in the case $F = \P^2$. Suppose for
simplicity that $\f$ is divisorial with an exceptional divisor $E$.
Using the results from the previous section (see (3.5)) we can consider
$\beta:\hat X\ra X$, the blow-up of $X$ along $F$ with the
exceptional divisor $\hat F$, $\alpha: \hat Z\ra Z$, the blow-up of
$Z$ at the maximal ideal of $z$ with the exceptional set $\hat G$,
and the map $\hat\f: \hat X\ra\hat Z$, which is a Fano-Mori contraction
supported by $-\hat F$. We note that, since $\hat
F=\P(N^*_F)$, we know the restriction of the map $\hat\f$ to
$\hat F$. Let $\hat E$ and $\hat S$ denote, respectively, the strict
transforms of $E$ and $S$.  All the objects are presented on the
following diagram.
$$\matrix
&&&\hat X\supset (\hat E,\hat F)&&&\cr
&&\beta\swarrow&&\searrow\hat\f&&\cr &&&&&&\cr &X\supset (E\supset
F)&&&&\hat Z\supset (\hat S,\hat G)&\cr &&&&&&\cr
&&\f\searrow&&\swarrow\alpha&&\cr
&&&Z\supset (S\ni z)&&&
\endmatrix
\leqno (*)$$

We provide a description of all these objects and, in
particular, of the singularities of $Z$ at $z$.  Actually, we can get
the description of $Z$ at $z$ by using merely (3.6) and (3.8), however, we
find it interesting to describe all the
objects which occur in this blow-up-blow-down construction.

The variety $\hat X$ is smooth and one can check (see (3.4)) that the map
$\hat\f_{|\hat F}$ is either the blow-up of a smooth three
dimensional quadric along a line or the blow-up of $\P^3$ along a
twisted cubic curve.  Consequently, the map $\hat\f$ is a divisorial
contraction with all non trivial fibers of dimension $1$. Therefore,
by (4.1), $\hat E$ is smooth and $\hat\f$ blows it down to a smooth
$\hat S$ in $\hat Z$ which is also smooth.

If $m$ is the multiplicity of $E$ along $F$, then we have the
numerical equivalence $$K_{\hat X}= \beta ^* K_X +\hat F = \beta^* E
+\hat F = \hat E + (m+1) \hat F$$ and also $K_{\hat X}= \hat\f^*
K_{\hat Z} + \hat E$.  Therefore $\hat\f^* K_{\hat Z} = (m+1) \hat F$
and, using the adjunction formula $(K_{\hat Z} +\hat G)_{|\hat G} =
K_{\hat G}$, we obtain that $(m+2)\hat G_{|\hat G} = \O(-3)$ or
$\O(-4)$ if $\hat G$ is a quadric or $\P^3$, respectively.  This
implies that the normal of $\hat G$ in $\hat Z$ is in both cases $\O
(-1)$ and also that $m = 1$, if $\hat G$ is a quadric, while $m = 2$
if $\hat G$ is $\P^3$.  Since the above multiplicity coincides with
the degree of the projection $\hat E\cap\hat F\ra F$ (which we know
because we know $\hat F=\P(N^*_{F/X})$) it follows that $\hat E$ and
$\hat F$ intersect transversally.  The description of the singularity
of $Z$ follows thus immediately (e.g.~one can apply (3.7) to
$\alpha$).
\par
Now we discuss the singularity of $S$ at $z$.  For this purpose we
consider the intersection curve $f:=\hat G\cap\hat S$. In both cases
$f\iso\P^1$ and the curve is either a line or a twisted cubic curve,
if respectively $\hat G\iso\Q^3$ or $\P^3$, and
$N_{f/\hat S}=\O(-1)$ or $\O(-3)$, if respectively
$G\iso\Q^3$ or $\P^3$. Moreover
$\alpha^{-1}(m_z)\cdot\O_{\hat S}=\O_{\hat S}(-f)$
so that $\alpha$ contracts $f$ to a normal singularity.
This provides the description of the singularity of $S$ at $z$.

\bigskip \noindent
(4.8). In the case when $F$ is smooth (i.e. if $F = \P^2$ or $\P^1 \times \P^1$)
one can easily construct examples by considering the total space
of the normal bundle, c.f. (3.7).
 More precisely let $(F,\Cal E)$ be one of the
following pairs:
$$\aligned (\P^2,\ (T\P^2(-1)\oplus\O(1))/\O),\ \ \ \  &(\P^2,\
\O^{\oplus 4}/\O(-1)^{\oplus2})\\
(\P^1\times \P^1,\ \O(1,0)\oplus\O(0,1)),\ \ \ &(\P^2,\ \O(1)\oplus\O(1)).
\endaligned$$
Let then $X := \P(\Cal E \oplus\O)$
and let $\xi$ denote the tautological line bundle on $X$.
The subvariety of $X$ determined by ${\Cal E}\oplus\O \ra \O \ra 0$
is isomorphic to $F$, we will call it again $F$; moreover
$N^*_{F/X} = {\Cal E}$.

Since $\Cal E\oplus \O$ is spanned by global sections
then, by definition, also $\xi$ is
spanned by global sections on $X$.
Let $\f: X \ra Z$ be the map associated to the complete
linear system $|\xi|$.
By the theorem of Leray and Hirsch one gets that
$\xi^4 = c_1^2({\Cal E}\oplus \O) - c_2({\Cal E}\oplus \O) > 0$.
This implies that the divisor $\xi$ is big so that
the map $\f$ is birational.

In the singular (also reducible) case it is more difficult to
construct examples. In \cite{A-W2} it was done by taking a codimension 2
complete intersection in a projective bundle over the smooth 4
dimensional quadric $\Q_4$. More precisely we considered a spinor
bundle ${\Cal S}$ over the quadric $\Q_4$ and the projective bundle
$M:=\P(\O\oplus{\Cal S}(1))\ra \Q_4$.  Again, let $\xi$ denote the
relative hyperplane bundle --- the evaluation map associated to the
system $|\xi|$ is onto $\P^4$ and contract a section $Q_0$ to a point
$v\in \P^4$.  The quadric $\Q^4$ parameterizes planes in $\P^4$
containing $v$ and $M$ is the incidence variety of points on the
planes.  We consider the linear system $\Lambda=|\xi+p^*\O(1)|$
over $M$ and $X$ was taken to be a complete intersection of two
divisors from $\Lambda$. The intersection $F:=X\cap Q_0$ is a linear
section of the 4-dimensional quadric. The codimension 2 linear section
of $Q_0\iso \Q^4$ can be either a smooth quadric
$\F_0=\P^1\times\P^1$ or a quadric cone $\S^2$, or a union of two
planes meeting along a line.  It can be shown that each of these cases can
be realized so that the complete intersection $X$ is smooth.
We refer the reader to \cite{A-W2} for a detailed
construction as well as for examples in the fiber type case.
\medskip

While closing this section let us note that a similar classification is
expected for contractions of fiber type from a smooth 4 dimensional
projective variety with an isolated two dimensional fiber. The list
of possible fibers is given in (4.3).  (As we already noted, similar results
were proved also by Kachi \cite{Kac}.) However, the understanding of the
normal of such a fiber requires a new approach. This is because the
restriction of $\f$ to a $X'\in|- K_X|$ is not a contraction but
generically a double covering.  In such a case the {\sl ascending
lemma} (3.11) has to be replaced by a {\sl trace argument} which, in
short, is the following. Since we can not lift sections from $X'$ to
$X$ because the restriction $\f_*\O_X\ra\f_*\O_{X'}$ is not
surjective any more, we consider the trace map $tr:
\f_*\O_{X'}\ra\O_Z=\f_*\O_X$. The map $tr$ has a nontrivial kernel
associated to the ramification of the (generically $2:1$) map $X'\ra
Z$.  It turns out that the kernel is not too large so that one can
lift up enough sections this way in order to relate properties of
fibers of $\f$ and $\f'$ and to prove the nefness of $N^*_{F/X}$ (see
\cite{A-W2}, (5.9.5)). Consequently, knowing the classification of nef
bundles over
$F$ with small Chern classes and using the results of Section 3 one can describe
the structure of $\f$ around $F$. In \cite{A-W2} we get the following:

\proclaim {Theorem (4.9)}
Let $\f:X\ra Z$ be a Fano-Mori contraction of a smooth 4-fold to a normal
3-fold.
Suppose that $F=\f^{-1}(z)$ is an irreducible isolated
2-dimensional fiber. If $F\iso \P^2$ then
$N^*_{F/X}$ is either $\O^3/\O(-2)$ or $T\P^2(-1)$. If $F$ is an
irreducible quadric
then $N^*_{F/X}$ is a pullback of $T\P^2(-1)$ via some double covering
of $\P^2$. In these cases the fiber structure $\tilde F$ coincides with
the geometric structure on $F$ and $Z$ is smooth at $z$.
\endproclaim

Let us note however that in the fiber case there is
also a number of possible exceptional fibers in (4.3) which are not locally
complete
intersection and they have to be dealt separately.  We expect to
discuss it in a forthcoming paper.

\bigskip

\subhead 5. Applications and extensions \endsubhead

In the present section we want to mention few applications of
techniques and results described previously.

Some results in this section are closely related to
{\sl adjunction theory of projective varieties},
a very classical theory (see \cite{C-E}), which was
reproposed and improved in modern time by A.~J.~Sommese and his
school.  One of the goals of this theory is to describe varieties $X$
polarized by an ample line bundle $L$ by means of the Fano-Mori
contraction supported by $K_X + rL$ where $r$ is the nef value of the
pair $(X,L)$.  If $X$ is smooth and $r \geq (n-2)$ then this goal is
achieved and we refer the reader to the new book \cite{B-S3} for an
overview of the theory; see also the next theorem (5.6). More generally
it is also achieved when the nef value is large with respect to the dimension
of fibers of $\f$.

The first result is a sequel of Theorem (1.10) in Section 1
and for $r = 1$ it reduces to Theorem (4.1).

\proclaim {Theorem (5.1)} Let $\f: X \ra Z$ be a Fano-Mori contraction
of a smooth variety $X$ and let $F=\f^{-1}(z)$ be a fiber.
Assume that $\f$ is supported by $K_X+rL$, with $L$ a $\f$-ample line
bundle on $X$.
\item {(1)} If $dim F \leq (r-1)$ then  $Z$ is smooth at $z$ and $\f$
is a projective bundle in a neighborhood of $F$.
\item{(2)} If $dimF=r$ then, after possible shrinking of $Z$ and restricting
$\f$ to a neighborhood of $F$, $Z$ is smooth and
\itemitem{(i)} if $\f$ is birational then $\f$ blows a
smooth divisor $E\supset X$
to a smooth codimension $r-1$ subvariety $S\supset Z$,
\itemitem{(ii)} if $\f$ is of fiber type and $dimZ=dimX-r$ then $\f$ is a
quadric bundle,
\itemitem{(iii)} if it is of fiber type and $dimZ = dimX-r+1$ then $r\leq
dimX/2$, $F =\P^r$ and the general fiber is $\P^{(r-1)}$.
\endproclaim

The proof of (5.1.2) i) and iii) is an application of
the relative base point freeness (1.11) together with slicing (1.14)
which, by induction, reduces to the result on good contractions
with 1-dimensional fibers (4.1). A detailed proof is contained in
the paper \cite{A-W1}. The Theorem (5.1.1) was proved by T. Fujita
(see \cite{Fu1},  (2.12)) also by a slicing argument; the proof of (5.5.2.ii)
was done in the paper \cite{A-B-W2} with the help of the proposition (1.3).

\medskip
The theorem (5.1) has many nice consequences, it opens further
generalizations and
conjectures. For instance an immediate consequence of the theorem (5.1.2) (i)
is the following generalization of the theorem 1.1 in \cite{E-S}.

\proclaim {Corollary (5.2)}
Let $\f: X \ra Y$ be a proper birational morphism between smooth varieties.
Put $S = \{y \in Y |\  {\hbox {\rm dim}}\f^{-1}(y) \geq 1 \}$, $E =
\f^{-1}(S)$.  Assume that
$(E)_{red}$ is an irreducible divisor and also that all the fibers of $\f_{|E}$
have the same dimension. Then $S$ is smooth and
$\f$ is the blow-up of $Y$ at $S$.
\endproclaim

Let us note that the case (iii) of the above theorem is related to the following
{\sl scroll conjecture} by Beltrametti and Sommese
\proclaim {Conjecture (5.3)}  Let $\f: X \ra Z$ be a Fano-Mori contraction
of a smooth variety $X$ supported by $K_X+rL$
and assume that $r = dimX - dim Z +1$ (the variety $X$ with such a map
$\f$ is then called a scroll or an adjunction scroll). If $dim X \geq
2dimZ -1$ then $Z$ is smooth and $\f$ is a projective bundle over
$Z$. If $dim X\geq 2dimZ-2$ then $Z$ is smooth.
\endproclaim
The first part of the conjecture is known if $dimX\geq 2dimZ$ and $L$ is very
ample by Ein \cite{Ei}, see also \cite{Wi3}.
It is immediate to see that the theorem implies the conjecture
if the special fibers have dimension not greater then $r$. In particular
the conjecture
is true if $dimZ \leq 2$ and also if $dimZ \leq 3$ and $\f$ is elementary.

\medskip \noindent
(5.4). Note that in the case (iii) of the theorem (5.3.2) we can define
a rank $r$ sheaf on $Z$ by setting $\E=\f_*(L)$.
The singular locus of the sheaf $\E$, call it $S$,
is the locus of $r$-dimensional fibers of $\f$.
Outside of $S$ the sheaf $\E$ is locally free of rank $r$
and at any point $z\in S$ we have the following local presentation of $\E$
$$0\ra \O_{Z,z}\ra\O_{Z,z}^{r+1}\ra\E_z\ra 0;$$
moreover $X =Proj_Z(Sym(\E))$.
Such sheaves are called {\sl B\v{a}nica sheaves}, they are discussed in
\cite{B-W}.

\smallskip
We suggest the following extension of Beltrametti-Sommese's conjecture:

\proclaim {Conjecture (5.4.1)} If $\f:X\ra Z$ is a smooth
scroll which is not a projective bundle and
$dimX=2dimZ-2$ then $X$ is the projectivization of a B\v{a}nica sheaf.
\endproclaim

\medskip

Similarly as in the case of (5.1), if we apply base point freeness
together with the slicing technique to results on birational
contractions of 4-folds, (4.7), then we get (see \cite{A-W2})

\proclaim {Proposition (5.5)}
Let $\f:X\ra Z$ be a Fano-Mori contraction of a smooth $n$-fold supported
by $K_X+(n-3)L$. Suppose that $\f$ is birational and that $F$ is an
isolated fiber of dimension $n-2$ (i.e. all neighboring fibers are
of dimension $\leq n-3$). If $n\geq 5$ then the contraction $\f$ is
small and $F$ is an isolated non-trivial fiber. More precisely
$F\iso\P^{n-2}$ and $N_{F/X}\iso\O(-1)\oplus\O(-1)$, and there exists
a flip of $\f$ (see e.g. \cite{A-B-W2} for the description of the flip).
\endproclaim

Possible extensions of the techniques discussed in this paper lead also
to the singular case. For example, using base-point-freeness one can
extend the results of Cutkovsky about contractions of 3-dimensional
varieties with terminal singularities --- see \cite{Cu} where also the
factoriality
was assumed. This was done in \cite{A1} and \cite{A2} and the following is one
of the results proved in these papers (see also a further development
in \cite{Me}); it is an extension
of Beltrametti-Sommese theory of first and second reduction to the case
of terminal singularities (see \cite{B-S2} and \cite{B-S3}).

\proclaim {Theorem (5.6)} Let $\f: X\ra Z$ be a Fano-Mori birational
contraction of
a projective variety with at most terminal singularities onto a
projective normal variety.
Let $L$ be a $\f$-ample line bundle such that $K_X +rL$ is a good supporting
divisor of $\f$.
\item{(1)} If $r > (n-2)$ then $r = (n-1)$ and $\f$ is the contraction
of a finite number of disjoint divisors $E_i$ such that $E_i \iso \P^{(n-1)}$,
each $E_i$ is contained in the smooth locus of $X$ and $N_{E_i/X} = \O(-1)$.
Moreover points $\f(E_i)$ are smooth in $Z$ and $Z$ has the same
singularities as $X$.
The contraction $\f$ is called the first reduction of the pair $(X,L)$.
\item{(2)} If $r =(n-2)$ assume also that $X$ is factorial
and $n >3$. Then $\f$ is the simultaneous
contraction of a finite number of disjoint divisors $D_i$ isomorphic to one
of the following:
\itemitem{(i)} $\P^{(n-1)}$ with normal bundle $\O(-2)$,
\itemitem{(ii)} $\Q^{(n-1)} \subset \P^n$ a hyperquadric (possibly singular)
with normal bundle $\O(-1)$,
\itemitem{(iii)} $D_i$ is the exceptional divisor of the blow-up of a curve $C$
contained as a complete intersection in the smooth locus of $Z$.
\item{} Note that new singularities appear on $Z$, namely the points $\f(D_i)$
where $D_i$ are divisors of type (i) and (ii).
The contraction $\f$ is called the second reduction of the pair $(X,L)$.
\endproclaim

One possible direction for generalization is to consider varieties
polarized by an ample vector bundle $\E$ whose rank is ``large''.  The
developments in this direction were stimulated initially by the
Mori's proof of Hartshorne conjecture (see \cite{Mo1}) and then by a
series of conjectures by Mukai.  Let us note that the present study
of the local structure of good contractions brings another
perspective to this problem. Namely, as we note in the previous
section, examples of Fano-Mori contractions can be built by using
projective bundles. Moreover, in the classification of conormal
bundles of contractions the knowledge of nef or spanned bundles
with small first Chern class played an essential role.

There is a number of papers devoted to this problem:
\cite{Wi1}, \cite{Fu3}, \cite{Pe}, \cite{Ye-Zh}, \cite{S-W1, 2},
\cite{Zh1}, \cite{A-B-W1}, \cite{A-M}. The following result was proved in
\cite{A-M}
and it covers most of the previous structural results.

\proclaim {Proposition (5.7)} Let $X$ be a smooth projective variety
and $\E$ be
an ample vector bundle of rank $r$ on $X$. Assume that $K_X+det\E$ is nef and
let $\pi:X\ra W$ be the
contraction supported by $K_X+det\E$.
\item{(1)} If $dim F\leq r-1$ for any fiber $F$ of $\pi$ and $r\geq (n+1)/2$
then $W$ is smooth, $\pi$ is of fiber type,
$F = \P^{r-1}$ and $\E_{|F}=\oplus^r\O(1)$.
\item{(2)}If $\pi$ is a divisorial elementary contraction,
with exceptional divisor $D$,
and  $dim F\leq r$ for all fibers $F$
then $W$ is smooth, $\pi$ is the blow up of a
smooth subvariety $B: = \pi (D)$ and
$\E =\pi^*\E^{\prime}\otimes[-D]$, for
some ample vector bundle $\E^{\prime}$ on $W$.
\endproclaim

The proof of this result is via Fano-Mori contractions of the variety
$\P(\E)$.  The proposition (5.7.1) holds also in the more general case
when $X$ has log terminal singularities and $codim Sing(X)>dim W$.
This was also proved in \cite{A-M}; the particular case $r= (n+1)$ was
proved independently by Zhang (see \cite{Zh1}).

\refstyle{A}

\enddocument